\documentclass[]{amsart}

\usepackage{graphicx}

\newcommand{\w}[1]{\ensuremath{\mathbf{\omega}_{#1}}}

\newtheorem{prop}{Proposition}

\newtheorem{fact}{Fact}
\newtheorem{conj}{Conjecture}

\theoremstyle{definition}
\newtheorem{defn}{Definition}

\newcommand{\p}[2]{\ensuremath{p^{#1}_{#2}}}
\newcommand{\q}[2]{\ensuremath{q^{#1}_{#2}}}

\newcommand{\MM}[2]{\ensuremath{\mathcal{M}^{#1}_{#2}}}
\newcommand{\LL}[2]{\ensuremath{\mathcal{L}^{#1}_{#2}}}

\newcommand{\CC}[1]{\ensuremath{\mathbf{C}^{#1}}}
\newcommand{\CP}[1]{\ensuremath{\mathbf{CP}^{#1}}}

\newcommand{\R}[1]{\ensuremath{\mathbf{R}^{#1}}}
\newcommand{\RP}[1]{\ensuremath{\mathbf{RP}^{#1}}}
\newcommand{\RR}{\ensuremath{\mathcal{R}}}
\newcommand{\HH}{\ensuremath{\mathcal{H}}}
\newcommand{\QQ}[2]{\ensuremath{\mathcal{Q}^{#1}_{#2}}}

\newcommand{\Q}[1]{\ensuremath{Q_{#1}}}

\newcommand{\Alt}[1]{\ensuremath{\mathcal{A}_{#1}}}
\newcommand{\Sym}[1]{\ensuremath{\mathcal{S}_{#1}}}
\newcommand{\D}[1]{\ensuremath{\mathcal{D}_{#1}}}
\newcommand{\G}[1]{\ensuremath{\mathcal{G}_{#1}}}
\newcommand{\Z}[1]{\ensuremath{\mathbf{Z}_{#1}}}

\newcommand{\obar}[1]{\ensuremath{\overline{#1}}}

\newcommand{\arrGap}{5pt}
\newcommand{\tabGap}{8pt}

\newcommand{\fig}[1]{Figure~\ref{#1}}
\newcommand{\figs}{Figures}

\newcommand{\app}[1]{\textsc{Appendix~\ref{#1}}}

\newcommand{\tab}[1]{Table~\ref{#1}}

\newcommand{\basinWidth}{9cm}

\begin{document}

\title[Solving the quintic by iteration]
 {Solving the quintic\\ by iteration in three dimensions}
\author{Scott Crass}
\address{Mathematics Institute\\ 
University of Warwick\\
Coventry  CV4 7AL  UK (1999-2000)}
\email{scrassw@maths.warwick.ac.uk}
\address{
Department of Mathematics\\
SUNY College at Buffalo\\
Buffalo, NY  14222}
\email{crasssw@buffalostate.edu}
%\keywords{complex dynamics, equivariant maps, quintic, symmetric group}
%\subjclass{51}

\begin{abstract}

The requirement for solving a polynomial is a means of breaking its
symmetry, which in the case of the quintic, is that of the symmetric
group \Sym{5}. Induced by its five-dimensional linear permutation
representation is a three-dimensional projective action. A mapping of
complex projective $3$-space with this \Sym{5} symmetry can provide the
requisite symmetry-breaking tool.

The article describes some of the \Sym{5} geometry in \CP{3} as well as
several maps with particularly elegant geometric and dynamical
properties. Using a rational map in degree six, it culminates with an
explicit algorithm for solving a general quintic.  In contrast to the
Doyle-McMullen procedure---three $1$-dimensional iterations, the present
solution employs one $3$-dimensional iteration.

\end{abstract}
\maketitle

%-------

\section{Overview}

In \cite{DM}, a solution to the quintic takes place in three iterative steps---a
\emph{tower} of algorithms each of which involves iteration in one complex
dimension.  Given almost any quintic $p$ and almost any initial point in \CC{}, 
the series of algorithms produces a root of $p$.  Their method is geometrically 
distinguished in that the tower has the \Sym{5} symmetry of the general quintic.  
Its central feature is a map on the Riemann sphere with icosahedral (\Alt{5}) 
symmetry.

The present paper describes a solution to a full measure's worth of quintics
that runs as a single iteration in three dimensions.  That the procedure
produces a root for almost any initial point in
complex projective $3$-space (\CP{3}) is conjectural at the
moment.  At its core is a map on \CP{3} with \Sym{5} symmetry.  Motivating this general 
project is a desire to develop solutions to equations that utilize geometrically
elegant dynamical systems.

The work unfolds in three stages: 1) some background geometry, 2) special maps
with \Sym{5} symmetry, and 3) a solution to the quintic based on the preceding
stages.

%----------

\textbf{Section \ref{sec:geom}: \Sym{5} geometry.}  
The setting here is \CP{3} 
upon which the symmetric group \Sym{5} acts. 
Finding a map with special \Sym{5} geometry requires some familiarity
with this action.  We will consider some features associated with the maps that
emerge in the second stage.  Indeed, the discovery of these maps derives from an 
awareness of the geometric landscape:

\begin{itemize}

\item coordinate systems

\item the structure of an \Sym{5}-invariant quadric surface 

\item the structure of certain special orbits of  points, lines, planes, and 
conics.  

\end{itemize}
In addition, the system of \Sym{5}-invariant polynomials plays a fundamental
role in the search for maps.

%----------

\textbf{Section \ref{sec:maps}: Maps with \Sym{5} symmetry.}
At this stage, we exploit our geometric understanding to discover empirically
several maps with special qualities.  Appearing here are families of maps 
associated with the icosahedron, the dodecahedron, and the complete graph on
five vertices. The known features of their geometric and dynamical behavior 
come under discussion.  However, they are not known to possess several desired 
properties.  In light of significant experimental evidence, I leave claims 
concerning these properties as conjectures.

%----------

\textbf{Section \ref{sec:quintic}: Dynamical solution to the quintic.}
Following the Doyle-McMullen framework, a special family of quintics corresponds 
to a \emph{rigid} family $\mathcal{E}$ of \Sym{5} maps on \CP{3}.  `Rigidity'
means that each member of $\mathcal{E}$ is conjugate to a single reference map
$f$ with elegant geometry and dynamics.  The solution is general since almost
any quintic $p$ transforms into the special family.  Thus, associated with $p$
is a map $g$ that we iterate.  Using \Sym{5} tools, its output---conjecturally,
a single \Sym{5} orbit---provides for an approximate solution to $\{p=0\}$.

A subsequent paper extends the method to the octic in a way that might
generalize to higher degree. \cite{CrassOct}

%-------

\section{\Sym{5} acts on \CP{3}} \label{sec:geom}

The permutation action of the \Sym{5} on \CC{5} preserves the hyperplane
$$\HH_x = \Biggl\{\sum_{k=1}^5 x_k =0 \Biggr\} \simeq \CC{4}$$ 
and, thereby, restricts to a faithful four-dimensional irreducible
representation.  (Since there will be two variables that describe the hyperplane,
the subscript $x$ appears here.)  This induces an \Sym{5} action on \CP{3}.   
Let \G{120} denote the corresponding subgroup of $\text{PGL}_4 \CC{}$.

%----------

\subsection{Coordinates}

For many purposes, the most perspicuous geometric description of \G{120}
employs five coordinates that sum to zero. One advantage is the simple
expression of the \G{120}-duality between points and planes. In general, for a
finite action \G{} whose matrix representatives are unitary, a point $a$ is 
\G{}-\emph{dual} to a hyperplane \LL{}{} if
$$\LL{}{} = \{ \obar{a} \cdot x = 0 \}.$$
Consequently, $a$ and
\LL{}{} have the same stabilizer in \G{}.  By the orthogonal
action of \Sym{5} on \CC{4}, a point
$$a = [a_1,a_2,a_3,a_4,a_5]_{\sum a_k = 0} \in \CP{3}$$ 
corresponds to the plane 
$$\{a \cdot x=0\}=\Biggl\{\sum_{k=1}^5 a_k\,x_k = 0\Biggr\}.$$
(Square brackets indicate a
point in projective space.)

A system of four $u$-coordinates also describes the hyperplane $\HH_u$. These
\emph{hyperplane coordinates} arise from the ``hermitian" change of
variable 
$$u = H\,x  \quad x = \obar{H^T} u \quad 
H=\frac{1}{\sqrt{5}}
\left( \begin{array}{ccccc} 1&\w{5}&\w{5}^2&\w{5}^3&\w{5}^4\\
1&\w{5}^2&\w{5}^4&\w{5}&\w{5}^3\\ 1&\w{5}^3&\w{5}&\w{5}^4&\w{5}^2\\
1&\w{5}^4&\w{5}^3&\w{5}^2&\w{5} \end{array} \right)
$$ 
where $\w{5}=e^{2\,\pi\,i/5}$ and the choice of scalar factor gives
\begin{equation} \label{eq:Hprod}
H\,\obar{H^T} = \begin{pmatrix} 
1&0&0&0\\
0&1&0&0\\
0&0&1&0\\
0&0&0&1 \end{pmatrix}
\quad
\obar{H^T} H = -\begin{pmatrix}
-4&1&1&1&1\\
1&-4&1&1&1\\
1&1&-4&1&1\\
1&1&1&-4&1\\
1&1&1&1&-4 \end{pmatrix}.
\end{equation}

%----------

\subsection{Invariant polynomials}

The fundamental result on symmetric functions states that the $n$
elementary symmetric functions of degrees one through $n$ generate the
ring of \Sym{n}-invariant polynomials. Since the \Sym{5} action on
\CP{3} occurs where the degree-$1$ symmetric polynomial vanishes, there
are four generating \G{120}-invariants. By Newton's identities, the
power sums
$$F_k(x) = \sum_{\ell=1}^5 x_{\ell}^k \qquad k=2, \dots, 5$$ 
also generate the \G{120} invariants. In hyperplane coordinates, these
are
\small
\begin{align*}
\Phi_2(u)=&\ F_2(\obar{H^T}u) = 2\,\bigl( u_1\,u_4 + u_2\,u_3 \bigr)\\
\Phi_3(u)=&\ \frac{3}{\sqrt{5}}\, \bigl( u_1\,u_2^2 + u_1^2\,u_3 +
u_3^2\,u_4 + u_2\,u_4^2 \bigr)\\[\arrGap]
\Phi_4(u)=&\ \frac{2}{5}\, \bigl( 2\,u_1^3\,u_2 + 3\,u_2^2\,u_3^2 +
2\,u_1\,u_3^3 + 2\,u_2^3\,u_4 + 12\,u_1\,u_2\,u_3\,u_4 + 3\,u_1^2\,u_4^2
+ 2\,u_3\,u_4^3 \bigr)\\[\arrGap]
\Phi_5(u)=&\ \frac{1}{5\,\sqrt{5}}\,\bigl(u_1^5 + u_2^5 +
20\,u_1\,u_2^3\,u_3 + 30\,u_1^2\,u_2\,u_3^2 + u_3^5 +
30\,u_1^2\,u_2^2\,u_4\\
&+ 20\,u_1^3\,u_3\,u_4 + 20\,u_2\,u_3^3\,u_4 + 30\,u_2^2\,u_3\,u_4^2 +
30\,u_1\,u_3^2\,u_4^2 + 20\,u_1\,u_2\,u_4^3 + u_4^5 \bigr).
\end{align*}
\normalsize

In classical invariant theory, relative invariants result from taking the
determinant of 1) the hessian $H(F)$ of an invariant $F$ and 2) the
``bordered hessian" $B(F,G)$ of two invariants $F$ and $G$
$$B(F,G) = 
\begin{pmatrix} 
&&&\frac{\partial{G}}{\partial{x_1}}\\
&H(F)&&\vdots\\
&&&\frac{\partial{G}}{\partial{x_n}}\\
\frac{\partial{G}}{\partial{x_1}}& \dots &
 \frac{\partial{G}}{\partial{x_n}}&0
\end{pmatrix}. $$
A polynomial $F$ is \emph{relatively invariant} if, for all $T \in \G{}$,
$$F \circ T = \alpha_T\,F$$
where $\alpha$ is a character on \G{}.

\begin{prop} \label{prop:hess}

Given $T \in \mathrm{GL}_n(\CC{})$ and invariants $F, G$,
\begin{align*}
|H(F(Tx))| =&\ |T|^{-2}\,|H(F(x))|\\
|B(F(Tx),G(Tx))| =&\ |T|^{-2}\,|B(F(x),G(x))|.
\end{align*}
Here, $|\cdot |$ indicates the determinant.

\end{prop} \noindent
For the permutation action of \Sym{5}, these give \emph{absolute}
invariants---the character is trivial. Thus, each is expressible in terms of the 
generators $\Phi_k$. The following result 
will serve a subsequent computational purpose.  (Note: Many of this work's 
results derive from calculation. For this purpose, I used \emph{Mathematica}.  
I will refer to them as ``Facts." )

\begin{fact} \label{fact:diffInv} 

With $G_4 = |H(\Phi_3)|$ and $G_5 = |B(\Phi_3,\Phi_2)|$, the ``power-sum"
invariants of degrees four and five are given by
\begin{align*}
\Phi_4 =&\ \frac{1}{324} \bigl(\Phi_2^2 - 5\,G_4 \bigr)\\
\Phi_5 =&\ \frac{1}{864} \bigl(720\,\Phi_2\,\Phi_3 + G_5 \bigr).
\end{align*}

\end{fact}

%----------

\subsection{Quadric surface}

The degree-$2$ invariant defines an \Sym{5}-invariant surface in \CP{3} 
$$\QQ{}{} =\{ \Phi_2 = 0 \}.$$ 
The quadratic form associated with \QQ{}{} is given by 
$$
\Phi_2 = 2\,\det U \quad \text{with} \quad 
U = \begin{pmatrix} u_1&-u_2\\u_3&u_4 \end{pmatrix}.
$$ 
Accordingly, \QQ{}{} is ruled by two families of lines 
$$ 
a^T\,U =
\begin{pmatrix} a_1&a_2 \end{pmatrix} 
 \begin{pmatrix} u_1&-u_2\\u_3&u_4 \end{pmatrix} 
  = \begin{pmatrix} 0&0 \end{pmatrix}
\quad U\,b = 
\begin{pmatrix} u_1&-u_2\\u_3&u_4\end{pmatrix} 
 \begin{pmatrix} b_1\\b_2 \end{pmatrix} 
  = \begin{pmatrix} 0\\0 \end{pmatrix}. 
$$
Alternatively, the ``$a$-ruling" is defined by 
$$U^T a = 0.$$
Each ruling forms a projective line $\CP{1}_a,\ \CP{1}_b$ respectively.         

Given a point $u=[u_1,u_2,u_3,u_4]$ on \QQ{}{}, the matrices $U$ and
$U^T$ each have rank one. Thus, distinct lines in $\CP{1}_a$ (or
$\CP{1}_b$) are skew while exactly one $a$-line and one $b$-line
intersect at $u$. This gives the quadric a $\CP{1}_a \times \CP{1}_b$
structure. (See \cite{HP}, Ch. XIII: Quadrics.)

Furthermore, as a set, each ruling has an \Alt{5} stabilizer \G{60} and,
hence, $\CP{1}_a$ and $\CP{1}_b$ have icosahedral geometry. The ``odd"
elements $\G{120} - \G{60}$ exchange the $a$-ruling with the
$b$-ruling.

%----------

\subsection{Special orbits}

The $3$-dimensional \Sym{5} action comes in both real and complex
versions. This means that, in the standard $x$ coordinates, \G{120} acts
on \RR---the \RP{3} of points with real components. \tab{tab:RP3} in
\app{app:tables} enumerates some special orbits contained in \RR\ while
\tab{tab:Q} describes elements of \QQ{}{} that are fixed by members of
\G{120}. For ease of expression, I will refer to special points (or
lines, planes, etc.) in terms of the orbit size: ``$20$-points"
($10$-lines, $5$-planes). Also, these points get a symbolic description
in reference to orbit size (superscript) and coordinate expression
(subscript).

Corresponding to each special point $a$ is the plane $\{a \cdot x=0\}$.
In the case of the 10-points 
$$[1,-1,0,0,0],\ \dots\ ,[0,0,0,1,-1],$$ 
there are 10-planes 
$$\{x_1 = x_2\},\ \dots\ ,\{x_4 = x_5\}$$ 
that are pointwise fixed by the involutions 
$$x_1 \leftrightarrow x_2,\ \dots\ , x_4\leftrightarrow x_5.$$
These ten transpositions generate \G{120} making it the projective image
of a \emph{real} or \emph{complex reflection group}. (See
\cite{ST}.)

Another noteworthy orbit is that of the five \Sym{4}-stable coordinate
planes
$$\LL{2}{5_i} = \{x_i = 0\} \quad i = 1, \dots, 5$$
as well as the five octahedral conics
$$\QQ{1}{i} = \QQ{}{} \cap \LL{2}{5_i}.$$
Some data for special two-dimensional orbits appear in \tab{tab:planes}.
I describe these sets in terms of dimension (superscript),
orbit-size (subscript), and coordinate expression (sub-subscript).

Finally, a number of special lines appear as intersections of the
$5$-planes and $10$-planes. \tab{tab:lines} summarizes the
situation.

%----------

\subsection{Configurations} \label{sec:config}

\begin{figure}[b]

\resizebox{2.4in}{!}{\includegraphics{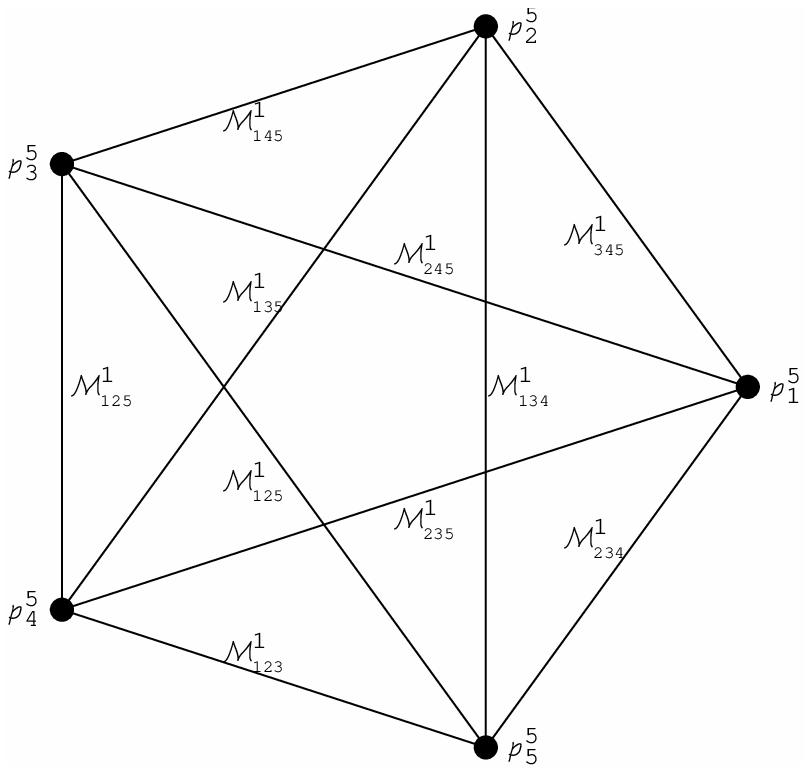}}
\resizebox{2.4in}{!}{\includegraphics{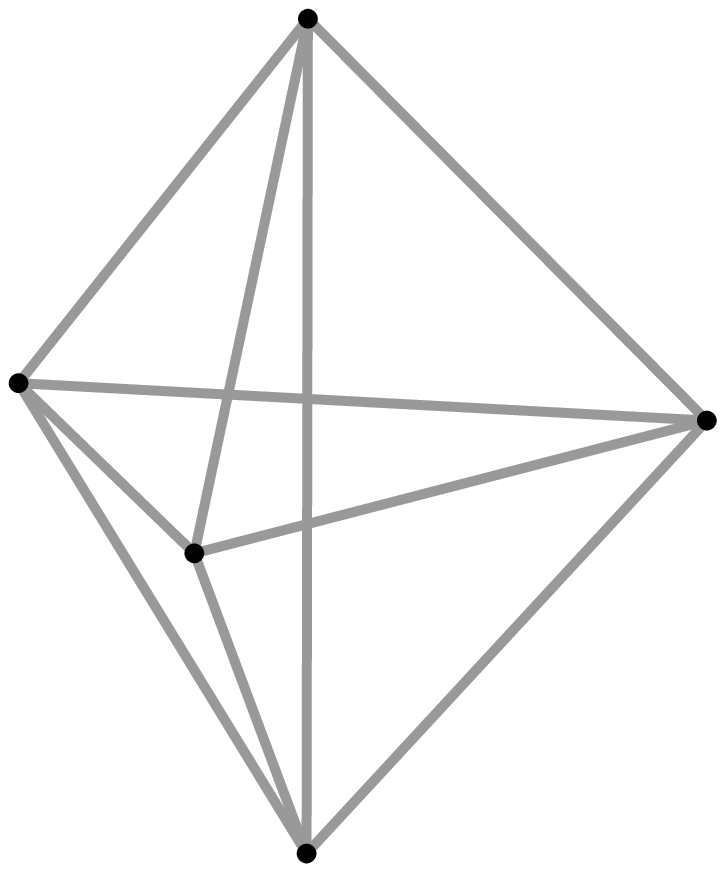}}

\caption{Configuration of $10$-lines and $5$-points}
 \label{fig:10lines}

\end{figure}

Some of the geometry that will have dynamical significance shows up in
various collections of lines. First, the 10-lines 
$$\MM{1}{10_{ijk}} = \LL{2}{10_{ij}} \cap \LL{2}{10_{ik}} \cap \LL{2}{10_{jk}}$$
form a complete graph on the 5-points. \fig{fig:10lines} illustrates
this in two ways. The pentagon-pentagram figure displays a $5$-fold
symmetry while the double pyramid exhibits the \D{3} structure of a
single $10$-line. (The
illustration suppresses the `10' subscript.)

The intersections of ``complementary" pairs of 10-planes yield an orbit
of 15-lines
$$
\LL{1}{15_{ij,k\ell}} = \LL{2}{10_{ij}} \cap  \LL{2}{10_{k\ell}} \qquad
\{i,j\} \cap \{k,\ell\} = \emptyset .
$$
This forms a graph on 15 vertices---the 5-points and 10-points
\p{10}{ij_2}.

\begin{itemize}

\item At a 5-point \p{5}{i}, there are three 15-lines
$$
\LL{1}{15_{jk,\ell m}}, \LL{1}{15_{j\ell,km}}, \LL{1}{15_{jm,k\ell}}
\quad i\neq j,k,\ell,m.
$$

\item On a 15-line \LL{1}{15_{jk,\ell m}}, there is one 5-point \p{5}{i}
where $i\neq j,k,\ell, m$.

\item At a 10-point \p{10}{ij_2}, there are three 15-lines 
$$
\LL{1}{15_{ij,k\ell}},\LL{1}{15_{ij,km}},\LL{1}{15_{ij,\ell m}}.
$$

\item On a 15-line \LL{1}{15_{ij,k\ell}} there are two 10-points
\p{10}{ij_2}, \p{10}{k\ell_2}.

\end{itemize}

Within each of the icosahedral rulings on \QQ{}{} there are three
special line-orbits. These correspond to the 12 vertices, 20
face-centers, and 30 edge-midpoints of the icosahedron. Intersections of
lines between rulings give special point structures.

\begin{figure}[b]
\resizebox{4in}{!}{
\begin{picture}(350,400)(-20,-250)
%\small
%
\thicklines
\put(0,-100){\line(1,1){200}}
\put(100,-200){\line(1,1){200}}
\put(0,-50){\line(1,0){300}}
\put(150,-200){\line(0,1){300}}
\thinlines
\put(0,0){\line(1,-1){200}}
\put(100,100){\line(1,-1){200}}
\put(150,50){\circle*{5}}
\put(150,-150){\circle*{5}}
\put(50,-50){\circle*{5}}
\put(250,-50){\circle*{5}}
\put(165,50){\q{20}{k\ell m_1}}
\put(165,-150){\q{20}{k\ell m_2}}
\put(152,90){\MM{1}{10_{k\ell m}}}
\put(45,-70){\q{20}{ij_1}}
\put(245,-70){\q{20}{ij_2}}
\put(0,-60){\LL{1}{10_{ij}}}
\put(130,-220){\vector(-1,1){20}}
\put(170,-220){\vector(1,1){20}}
\put(125,-230){$\QQ{}{}\ \cap\ T_{\q{20}{k\ell m_2}}\!\QQ{}{}$}
\put(130,115){\vector(-1,-1){20}}
\put(170,115){\vector(1,-1){20}}
\put(125,120){$\QQ{}{}\ \cap\ T_{\q{20}{k\ell m_1}}\!\QQ{}{}$}
\put(-15,-60){\vector(1,-2){15}}
\put(-15,-40){\vector(1,2){15}}
\put(-65,-50){$\QQ{}{}\ \cap\ T_{\q{20}{ij_1}}\!\QQ{}{}$}
\put(315,-60){\vector(-1,-2){15}}
\put(315,-40){\vector(-1,2){15}}
\put(305,-50){$\QQ{}{}\ \cap\ T_{\q{20}{ij_2}}\!\QQ{}{}$}
\end{picture}
}

\caption{Configuration of $40$-lines and $20$-points on \QQ{}{}.  At a 
$20$-point \q{20}{ij_-} or \q{20}{ijk_-} there are two
$40$-lines---one in each ruling on the quadric. This pair of lines is
the intersection of \QQ{}{} with the tangent plane to \QQ{}{} at the
respective $20$-point. Also indicated are the $10$-lines determined by a
pair of antipodal $20$-points.}

\label{fig:20pts}

\end{figure}
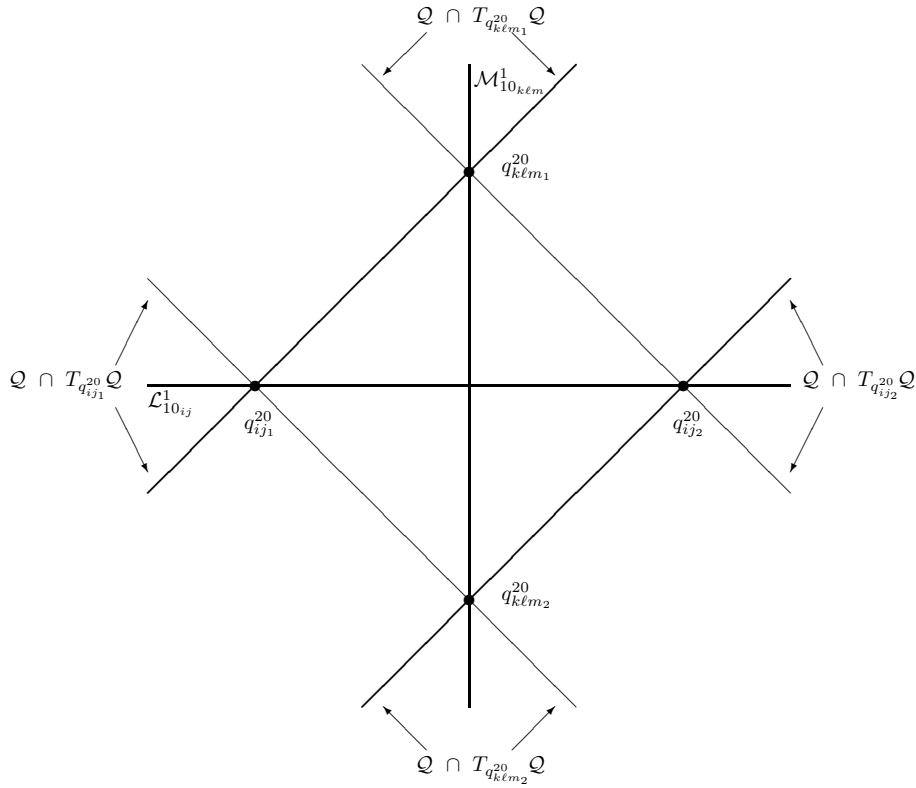

\begin{itemize}

\item Two $20$-line \G{60}-orbits form ten ``quadrilaterals" at two pairs
of $20$-points. (See \fig{fig:20pts}.)

\item Two $12$-line \G{60}-orbits form six quadrilaterals at
$24$-points.

\item Two $30$-line \G{60}-orbits orbits form $15$ quadrilaterals at two
pairs of $30$-points.

\end{itemize}
Since $\G{120}-\G{60}$, exchanges the orbits in $\CP{1}_a$ with those in 
$\CP{1}_b$, these give overall line-orbits of sizes $40$, $24$, and $60$.
%%%%%%%%%%%%%%
%%%%%%%%%%%%%%
%%%%%%%%%%%%%%

%-------

\section{Equivariant maps} \label{sec:maps}

The primary tool to be used in solving the general quintic is a rational
map 
$$f:\CP{3} \longrightarrow \CP{3}$$
with \Sym{5} symmetry. In algebraic terms, this means that 
$$f \circ T = T\circ f \quad \text{for all}\ T \in \G{120}.$$
Furthermore, such an \emph{equivariant} map (or simply \emph{equivariant}) 
should have \emph{reliable dynamics}: its attractor

\begin{enumerate}

\item[1)] is a \emph{single} \G{120} orbit

\item[2)] has a corresponding basin with full measure in \CP{3}

\item[$2'$)] alternatively, has a corresponding basin that is dense in
\CP{3}.

\end{enumerate}

Recall that a point $a$ in a space $X$ is \emph{attracting} when, for all $x$ in some
neighborhood of $a$,
$$f^k(x) \longrightarrow a.$$
A point $s$ is \emph{superattracting in a direction} $L$ if the derivative 
$f'(s)$ has
a zero eigenvalue in the $L$ direction.
The \emph{basin of attraction} $B_a$ of $a$ is the set of all points attracted to
$a$;
$$B_a = \{x \in X : f^k(x) \longrightarrow a\}.$$
Also, the \emph{attractor} of $f$ is the set of all attracting points. 

%----------

\subsection{Basic maps}

A finite group action \G{} on \CC{n} induces an action on the associated
exterior algebra. Moreover, \G{}-invariant $(n-1)$-forms
correspond to \G{}-equivariant maps. \cite{Crass}  Briefly,
let
$$dZ^I = (-1)^{\sigma_I} dz_{i_1} \wedge \dots \wedge dz_{i_{n-1}}$$
where $I$ is the ordered set 
$$\{i_1, \dots, i_{n-1}\}\qquad i_1 < \dots < i_{n-1},$$
$\hat{I}$ is the single index in $\{1, \dots, n\} - I,$ and $\sigma_I$ is the 
sign of the permutation
$$
\begin{pmatrix}
1&2& \dots& n\\
\hat{I}&i_1& \dots & i_{n-1}
\end{pmatrix}.
$$
If
$$\phi(z) = \sum_{\hat{I}=1}^n f_{\hat{I}}(z) dz^I$$
is a \G{}-invariant $(n-1)$-form, then the map
$$f(z) = (f_1(z), \dots, f_n(z))$$
is \emph{relatively} \G{}-equivariant (a multiplicative character appears under
the action of \G{} on $f$).

For a reflection group, the number of generating $0$-forms (i.e.,
polynomials) is the dimension of the action. \cite[p. 282]{ST}
From a result in complex reflection groups, this is also the number of
generating $1$-forms and $(n-1)$-forms. \cite[p. 232]{OT}
Indeed, the $1$-forms are exterior derivatives of the $0$-forms while
the $(n-1)$-forms are wedge products of $1$-forms.

\begin{prop} \label{prop:genMaps}

With $X^k_i = - 4\,x_i^k + \sum_{j\neq i} x_j^k $, the four maps 
$$
f_k(x) = \bigl[X^k_1,X^k_2,X^k_3,X^k_4,X^k_5\bigr] 
\qquad k=1, \dots, 4
$$
generate the module of \G{120} equivariants over the ring of 
\G{120}-invariants.

\end{prop} \noindent

These maps are projections onto the hyperplane $\HH_x$ along
$[1,1,1,1,1]$ of the power maps
$$\bigl[x_1^k,x_2^k,x_3^k,x_4^k,x_5^k\bigr].$$

\begin{prop}

Under an orthogonal action an invariant $F(x)$ gives rise to an
equivariant $f(x)$ by means of a formal gradient
$$f(x) = \nabla_x F(x) = 
\biggl[\frac{\partial F}{\partial x_1}(x),\dots, 
 \frac{\partial F}{\partial x_n}(x)\biggr].$$

\end{prop}

\begin{proof}

For a homogeneous polynomial $F(x)$ of degree $m$, the Euler identity
gives
$$
m\,F(x) = \nabla_x F(x)^T x = \nabla_x F(x) \cdot x = f(x)\cdot x.
$$
Invariance of $F$ yields
$$m\,F(x) = m\,F(Ax) = \nabla_x F(Ax)^T Ax.$$
Using an auxiliary variable $y$,
$$\nabla_x F(Ax) = A^T\,\nabla_y F(y)|_{y = A\,x} =  A^T\,f(y) = A^T\,f(A\,x).$$
By orthogonality of $A$,
$$m\,F(x) = A^{-1}\,f(A\,x)\cdot x.$$
Equating expressions for $m\,F(x)$ reveals equivariance:
$$A^{-1} f(Ax) = f(x).$$
\end{proof}
Note that the \Sym{5}-equivariant $f_k(x)$ is \emph{not} equal to $\nabla_x
F_{k+1}(x)$, but is a multiple of 
$$\nabla F_{k+1}(x)|_{x_i^k =X_i^k}.$$
While this may be a source of confusion, it does not cause problems, since we
are working on the hyperplane $\HH_x$.  When using hyperplane coordinates on
$\HH_u$, the discrepancy disappears.

A map on $\HH_x$ produces 
$$\phi(u) = H f(\obar{H^T} u)$$
on $\HH_u$. Expressing the generating $u$-equivariants 
$$\phi_k(u) = H f_k(\obar{H^T} u)$$
in terms of the basic $u$-invariants $\Phi_k(u)$ will be
useful.

\begin{defn}

Let 
$$R = \begin{pmatrix} 
0&0&0&1\\ 
0&0&1&0\\
0&1&0&0\\
1&0&0&0
\end{pmatrix} \quad
\text{and} \quad
\nabla_u^r F(u) = R\,\nabla_u F(u)$$
represent the \emph{reversed identity} and \emph{reversed gradient}.

\end{defn}

\begin{prop} \label{prop:Hmaps}

In $\HH_u$ coordinates, the map $\phi(u) = H\,f(\obar{H^T} u)$
is given by 
$$\phi(u) = \nabla_u^r \Phi(u)$$
where $\Phi(u) = F(\obar{H^T} u) = F(x)$.

\end{prop}

\begin{proof}

For the change of variable $u = H\,x$ and $x = \obar{H^T}\,u$, the chain
rule yields
$$f(x) = \nabla_x F(x) = \nabla_x \Phi(u) = H^T\,\nabla_u \Phi(u).$$
Since $H H^T = R$,
\begin{align*}
H\,f(x) =&\ H\,H^T\,\nabla_u \Phi(u)\\
H\,f(\obar{H^T}\,u) =&\ R\,\nabla_u \Phi(u)\\
\phi(u) =&\ \nabla_u^r \Phi(u).
\end{align*}

\end{proof}

Thus, the basic maps in $u$ are
$$\phi_k(u) = -\frac{5}{k+1}\,\nabla_u^r \Phi_{k+1}(u).$$
Explicitly,
\small
\begin{align*} 
\phi_1(u) =&\ 2\,\bigl[u_1,u_2,u_3,u_4 \bigr]\\
\phi_2(u) =&\
\frac{3}{\sqrt{5}}\,\bigl[
u_3^2 + 2\,u_2\,u_4, u_1^2 + 2\,u_3\,u_4, 
2\,u_1\,u_2 + u_4^2, u_2^2 + 2\,u_1\,u_3 
\bigr]\\
\phi_3(u) =&\
\frac{4}{5}\,\bigl[
u_2^3 + 6\,u_1\,u_2\,u_3 + 3\,u_1^2\,u_4 + 3\,u_3\,u_4^2,
3\,u_2^2\,u_3 + 3\,u_1\,u_3^2 + 6\,u_1\,u_2\,u_4 + u_4^3,\\
&u_1^3 + 3\,u_2\,u_3^2 + 3\,u_2^2\,u_4 + 6\,u_1\,u_3\,u_4,
3\,u_1^2\,u_2 + u_3^3 + 6\,u_2\,u_3\,u_4 + 3\,u_1\,u_4^2
\bigr]\\
\phi_4(u) =&\
\frac{6}{\sqrt{5}}\,\bigl[
4\,u_1^2\,u_2^2 + 4\,u_1^3\,u_3 + 4\,u_2\,u_3^3 + 12\,u_2^2\,u_3\,u_4 +
12\,u_1\,u_3^2\,u_4 + 12\,u_1\,u_2\,u_4^2 + u_4^4,\\
&4\,u_1\,u_2^3 + 12\,u_1^2\,u_2\,u_3 + u_3^4 + 4\,u_1^3\,u_4 +
12\,u_2\,u_3^2\,u_4 + 6\,u_2^2\,u_4^2 + 12\,u_1\,u_3\,u_4^2,\\
&u_2^4 + 12\,u_1\,u_2^2\,u_3 + 6\,u_1^2\,u_3^2 + 12\,u_1^2\,u_2\,u_4 +
4\,u_3^3\,u_4 + 12\,u_2\,u_3\,u_4^2 + 4\,u_1\,u_4^3,\\
&u_1^4 + 4\,u_2^3\,u_3 + 12\,u_1\,u_2\,u_3^2 + 12\,u_1\,u_2^2\,u_4 +
12\,u_1^2\,u_3\,u_4 + 6\,u_3^2\,u_4^2 + 4\,u_2\,u_4^3 \bigr].
\end{align*}
\normalsize

%----------

\subsection{A fixed point property}

For a \G{120}-equivariant $f$ and a point $a$ that an element $T \in
\G{120}$ fixes, 
$$T\,f(a) = f(Ta) = f(a).$$
Hence, equivariants preserve fixed points of a group element. 

Being pointwise fixed by the involution
$$x_i \longleftrightarrow x_j,$$
a $10$-plane 
$$\LL{2}{10_{ij}} = \{x_i - x_j =0\}$$
either maps to itself or collapses to its companion $10$-point
$$
\p{10}{ij_1} = 
[\dots 0 \dots,\overbrace{1}^{i},\dots 0 \dots,
 \overbrace{-1}^{j},\dots 0 \dots] 
\notin \LL{2}{10_{ij}}.
$$
In the former generic case, the map preserves the $10$-line
and $15$-line orbits \MM{1}{10_{ij}} and \LL{1}{15_{ij,k\ell}} that are
intersections of $10$-planes.

%----------

\subsection{Families of equivariants}

The \G{120} equivariants form a module over the \G{120} invariants for which
degree provides a grading. This means that for an invariant $F_\ell$ and
equivariant $g_m$ of degrees $\ell$ and $m$, the product 
$$F_\ell \cdot g_m$$ 
is an equivariant of degree $\ell + m$. When looking for a map in a certain
degree $k$ with special geometric or dynamical properties, my approach is to
express the entire family of ``$k$-maps" and by manipulation of parameters,
locate a subfamily with the desired behavior.

%----------

\subsection{Quadric-preserving maps}

The rich geometry of the quadric \QQ{}{} provides an intriguing setting
for dynamical exploration. Are there \Sym{5}-symmetric maps that send
\QQ{}{} to itself? If so, how do they behave on and off \QQ{}{}? I will
describe discoveries of two species of such maps: one associated with
the icosahedron and the other with the octahedron.

%-------------

\subsubsection*{Maps that preserve icosahedral rulings}

Were a \G{120}-equivariant to preserve the \Alt{5} rulings on \QQ{}{},
its restriction to either ruling $\CP{1}_a$ or $\CP{1}_b$ would express
itself in terms of the basic equivariants under the one-dimensional
icosahedral action.  Such maps occur in
degrees $11$, $19$, and $29$. \cite[p. 166]{DM} Consequently, the 
$20$-parameter family of $11$-maps comes under scrutiny:
\small
\begin{align*} 
f_{11} =\ 
&(\alpha_{1}\,F_2^5 + \alpha_{2} F_2^2\,F_3^2 + \alpha_{3}\,F_2^3\,F_4 +
\alpha_{4}\,F_3^2\,F_4 + \alpha_{5}\,F_2\,F_4^2 +
\alpha_{6}\,F_2\,F_3\,F_5 + \alpha_{7}\,F_5^2)\,f_1 \\
&+ (\alpha_{8}\,F_2^3\,F_3 + \alpha_{9}\,F_3^3 +
\alpha_{10}\,F_2\,F_3\,F_4 + \alpha_{11}\,F_2^2\,F_5 +
\alpha_{12}\,F_4\,F_5) f_2 \\ 
&+ (\alpha_{13}\,F_2^4 + \alpha_{14}\,F_2\,F_3^2 +
\alpha_{15}\,F_2^2\,F_4 + \alpha_{16}\,F_4^2 +
\alpha_{17}\,F_3\,F_5)\,f_3 \\
&+ (\alpha_{18}\,F_2^2\,F_3 + \alpha_{19}\,F_3\,F_4 +
\alpha_{20}\,F_2\,F_5)\,f_4. 
\end{align*}
\normalsize

From the geometric description of the
icosahedral $11$-map on $\CP{1}_a$ or $\CP{1}_b$ (\cite[p. 163]{DM}), a 
ruling-preserving
$11$-map would exchange antipodal pairs of $20$-lines
$\{\LL{a}{20_1},\LL{a}{20_2}\}$ or $\{\LL{b}{20_1},\LL{b}{20_2}\}$ and 
$30$-lines while fixing $12$-lines. (Recall that these are \G{60} orbits.)
Imposed on the configurations described in Section~\ref{sec:config},
these conditions require analogous behavior at the associated points:
\begin{align*} 
\q{20}{ij_1} &\longleftrightarrow \q{20}{ij_2}& 
\q{20}{ijk_1}&\longleftrightarrow \q{20}{ijk_2}\\ 
\q{24}{ijk\ell} &\longrightarrow \q{24}{ijk\ell}\\ 
\q{30}{i,jk_1} &\longleftrightarrow \q{30}{i,jk_2}&
\q{30}{ij,k\ell_1}&\longleftrightarrow \q{30}{ij,k\ell_2}. 
\end{align*}
The specified action occurs automatically for $\q{20}{ij_-},\
\q{24}{ijk\ell},\ \text{and}\ \q{30}{i,jk_-}.$ After solving two
linear equations associated with the remaining two conditions
$$
f_{11}(\q{20}{ijk_1}) =\ \q{20}{ijk_2}\qquad
f_{11}(\q{30}{ij,k\ell_1}) =\ \q{30}{ij,k\ell_2}
$$
as well as four linear equations 
\begin{equation*}
f_{11}(\LL{a}{20_1}) =\ \LL{a}{20_2}
\end{equation*}
that arrange for the exchange of an antipodal pair of $20$-lines in
either ruling we obtain a $13$-parameter family of ruling-preserving
maps
\small \begin{align*} 
g_{11} =&\ 
4\,\bigl( 16\,\alpha_1\,F_2^5 + 16\,\alpha_2\,F_2^2\,F_3^2 +
16\,\alpha_3\,F_2^3\,F_4 + 67\,F_3^2\,F_4 \\
&+\ 16\,\alpha_5\,F_2\,F_4^2 + 16\,\alpha_6\,F_2\,F_3\,F_5 
+ 45\,F_5^2 \bigr)\,f_1\\
&+\ 4\,\bigl( 16\,\alpha_8\,F_2^3\,F_3 + 16\,F_3^3 
+ 16\,\alpha_{10}\,F_2\,F_3\,F_4 + 16\,\alpha_{11}\,F_2^2\,F_5 
- 135\,F_4\,F_5 \bigr)\,f_2 \\
&+\ \bigl(64\,\alpha_{13}\,F_2^4 + 64\,\alpha_{14}\,F_2\,F_3^2 + 
64\,\alpha_{15}\,F_2^2\,F_4 + 405\,F_4^2 - 720\,F_3\,F_5 
\bigr)\,f_3\\
&+\ 4\,\bigl(16\,\alpha_{18}\, F_2^2\,F_3 - 225\,F_3\,F_4 
+ 16\,\alpha_{20}\,F_2\,F_5 \bigr)\,f_4. 
\end{align*} \normalsize
When restricted to the ruling $\CP{1}_a$ and expressed in the
homogeneous \emph{ruling coordinates} $[a_1,a_2]$, the map has the elegant
appearance
\begin{equation*}
g_{11}|_{\CP{1}_a}: [a_1,a_2] \longrightarrow
\bigl[a_1\,\bigl(-a_1^{10} + 66\,a_1^5\,a_2^5 + 11\,a_2^{10} \bigr),
a_2\,\bigl(11\,a_1^{10} - 66\,a_1^5\,a_2^5 - a_2^{10} \bigr) \bigr].
\end{equation*}
Of course, the same form appears for the $b$-ruling.

Restricted to a ruling, the dynamics of each $g_{11}$ is completely
understood. The $20$-lines  are period-$2$ and the only elements of the
critical set. (Recall that $20$-lines in \QQ{}{} are dodecahedral vertices in
$\CP{1}_a$ or $\CP{1}_b$.)  This implies that almost every line in the ruling 
belongs to the basin of one of the ten pairs of the superattracting
set. (See \cite[p. pp. 166-167]{DM} and \fig{fig:dodec11} in
\app{app:basins}.) Thus, for almost every point $q_0$ on \QQ{}{}, there
is an``antipodal" pair of intersections \emph{between} $20$-lines in
each ruling toward which $g_{11}$ attracts the trajectory: 
$$
g^n_{11}(q_0) \longrightarrow \{ \LL{a}{20_1} \cap
\LL{b}{20_i},\LL{a}{20_2} \cap \LL{b}{20_j} \} \quad \{i,j\} = \{1,2\}.
$$
As a result, the global behavior of each $g_{11}$ depends on its
dynamics off \QQ{}{}. Should the quadric attract or repel? If \Q{}{}
were attracting, then the $400$ intersections of $20$-lines would
attract in all directions. One way to arrange for this is to force these
points to be critical in the off-quadric direction. However, this
situation does not conform to the model of \emph{reliable} dynamics. The
attractor would not be a single \G{120} orbit of points, though it might be the
set of intersections of a single line-orbit. I have not explored the
case of a repelling quadric. 

Interestingly, the quadric resists criticality. Computation reveals that
no member of $g_{11}$ is critical on all of \QQ{}{}.  Is there a geometric
reason for this?  The next example
reveals that this is not a universal trait of quadric-preserving maps. 

%-------------

\subsubsection*{An octahedral map}

Since the orbit of five planes \LL{2}{5_k} has fundamental geometric
significance, a map that preserves these sets might exhibit interesting
dynamics.  Arranging for this spends four of the twenty parameters of
the family $f_{11}$.

The intersection of a $5$-plane \LL{2}{5_k}and \QQ{}{} is a conic
\QQ{1}{k} with \Sym{4} symmetry and, thereby, octahedral structure. One
of the special equivariants for the octahedral action on \CP{1} is a
$5$-map that attracts almost every point to the eight face-centers---
vertices of the dual cube. Geometrically, the map stretches each face
$\mathsf{F}$ of the cube symmetrically over the five faces in the
complement of the face antipodal to $\mathsf{F}$. As a face stretches,
it makes a half-turn so that the vertices land on their antipodes. This
makes each vertex critical and period-$2$; locally, the map is squaring.
Since these are the only critical points, their basins have full
measure. (See \cite[p. 156]{DM} and \fig{fig:oct5} in
\app{app:basins}.) Under \G{120}, antipodal pairs of octahedral
face-centers are the $20$-points \q{20}{ij_1}, \q{20}{ij_2}.

The idea is to look for a reliable map with the $20$-points as its only
attractor. In degree five there are too few parameters for the purpose.
However, the $11$-maps provide enough freedom to arrange for elegant
geometry. The goal demands that the desired map $h_{11}$ preserve the
$5$-conics \QQ{1}{k} and then decay to the octahedral $5$-map there. One
way to realize this is to self-map the quadric \QQ{}{}. This takes six
of the remaining $16$ parameters the expenditure of one of which assures
that the $20$-points do not blow up.

Intriguingly, when any member $g_{11}$ of the resulting $10$-parameter
family restricts to \QQ{}{}, it decays into a $5$-map
$$ 
g_{11}|_{\QQ{}{}} =
-\frac{1}{2}\,F_3^2\,\bigl(2\,F_3\,f_2 - F_4\,f_1 \bigr)|_{\QQ{}{}}.
$$
This decadence occurs unexpectedly, since most octahedral $11$-maps exchange
pairs of face-centers and are non-degenerate.
When restricted to an ``affine" part of the quadric
$$\QQ{}{} \cap \{u_1 \neq 0\},$$
the maps have the simple form
$$
g_{11}|_{\QQ{}{} \cap \{u_1 \neq 0\}}:
(x,y) \longrightarrow 
\biggl( 
\frac{x^2+ 3\,y - 2\,x\,y^3}{2\,x + 3\,x^2\,y^2 - y^3},
\frac{3\,x^2 + 2\,y + x^3\,y^2}{1 + 2\,x^3\,y - 3\,x\,y^2}
\biggr).
$$
Is there a geometric
description of the restricted map?

Every member of the $g_{11}$ family preserves the \Sym{3}-symmetric conic
$$\QQ{}{} \cap \LL{2}{10_{ij}}$$ 
each of which contains a pair of $20$-points \q{20}{ij_1}, \q{20}{ij_2}.
In coordinates where these points are $0$ and $\infty$, 
$$ 
g_{11}|_{\QQ{}{} \cap \LL{2}{10_{ij}}}:
z \longrightarrow 
\frac{7\,\sqrt{5}\,z^3 + 5\,i}{z^2\,(5\,i\,z^3 + 7\,\sqrt{5})}. 
$$
Of course, the period-$2$ points $0$ and $\infty$ are critical. By
experiment, the remaining six critical points belong to their
superattracting basin. Such circumstances force almost every point on a
conic to belong to the basin.

Octahedral $11$-maps generically exchange antipodal pairs of vertices.
Such a pair corresponds to the $30$-points \q{30}{i,jk_1},
\q{30}{i,jk_2}. As a degenerate member of the family, the $5$-map fixes
these points. These conditions require each $g_{11}$ to blow up the
$30$-points. Also blowing up are the $24$-points.

Now the issue is behavior off \QQ{}{}. Since the desired attractor lies
on \QQ{}{} and the dynamics there appears to be reliable, a map for
which the quadric is itself attracting comes to mind. Because octahedral
face-centers are superattracting on the respective conics, each $g_{11}$
is critical at the associated $20$-points. The maps are also critical at
the blown-up points. Arranging for critical behavior at the three
quadric orbits consisting of the non-octahedral $20$-points, $30$-points
and the octahedral $60$-points costs three parameters. The result is a
seven parameter family of $11$-maps for which the entire quadric is
critical and each octahedral $20$-point is superattracting in three
directions.

Each of the $10$-lines \LL{1}{10_{ij}} contains a pair of antipodal
$20$-points. A map that 

\begin{itemize}

\item[1)] preserves these lines, 

\item[2)] attracts almost every point on the line to the $20$-points,
and 

\item[3)] is critical in the directions \emph{off} the line 

\end{itemize}
would act as a ``superattracting pipe" to the quadric. Expenditure of
four of the remaining seven parameters purchases a map with these
properties. Indeed, when restricted to each \LL{1}{10_{ij}}, the map is
$$z \longrightarrow -\frac{1}{z^2}$$
with the pair of $20$-points at $0$ and $\infty$.

The final three parameters allow for a map $h_{11}$ with a
non-attracting pipe to \QQ{}{} at the $10$-lines \MM{1}{10_{ijk}}:
\small \begin{align*} 
h_{11} =&\ 
\bigl(-21\,F_2^5 + 56\,F_2^2\,F_3^2 + 66\,F_2^3\,F_4 + 48\,F_3^2\,F_4 -
48\,F_2\,F_4^2 - 96\,F_2\,F_3\,F_5 \bigr)\,f_1\\ 
&-\ 24\,\bigl(4\,F_3^3 - 9\,F_2\,F_3\,F_4 + 3\,F_2^2\,F_5
\bigr)\,{f_2}
+\ 12\,\bigl(5\,F_2^4 + 8\,F_2\,F_3^2 - 10\,F_2^2\,F_4 \bigr)\,f_3\\
&-\ 96\,F_2^2\,F_3\,f_4.
\end{align*} \normalsize
Such a line contains the pairs of $20$-points \q{20}{ijk_1},
\q{20}{ijk_2}. In coordinates where these points are $0$ and $\infty$,
the restriction of $h_{11}$ to \MM{1}{10_{ijk}} is
$$z \longrightarrow -\frac{1}{z^2}.$$
On \QQ{}{} these $20$-points are repelling.  Indeed, they belong to the
conics
$$
\QQ{}{} \cap \LL{2}{10_{ij}},\ \QQ{}{} \cap \LL{2}{10_{ik}},\ 
\QQ{}{} \cap \LL{2}{10_{jk}}
$$
on which the basins of the pair of $20$-points \q{20}{ij_-},
\q{20}{ik_-}, \q{20}{jk_-} have full measure. Experiment reveals that
nearby points belong to the basins of the other $20$-point orbit. 

Due to its geometry, $h_{11}$ preserves the various \CP{1} intersections
of $5$-planes and $10$-planes. The two such lines not yet considered are
the $15$-lines \LL{1}{15_{ij,k\ell}} and the $30$-lines
\LL{1}{30_{i,jk}}. In ``symmetrical" coordinates where the intersections
with \QQ{}{} are at $0$ and $\infty$, the restricted maps are:
\begin{align*}
&h_{11}|_{\LL{1}{15_{ij,k\ell}}}: 
z \longrightarrow \frac{19\,z^2 - 9}{z^2(9\,z^2 - 19)}\\
&h_{11}|_{\LL{1}{30_{i,jk}}}: 
z \longrightarrow -\frac{11\,z^2 + 9}{z^2(9\,z^2 + 11)}.
\end{align*}
In the former case, the map has attracting fixed points at the pair of
$10$-points \p{10}{ij_2}, \p{10}{k\ell_2} and a period-$2$
superattractor at \q{30}{ij,k\ell_1}, \q{30}{ij,k\ell_2}. Overall, these
are saddle points where the map repels off the line. A similar state of
affairs occurs on the $30$-lines. Here, the pair of attracting fixed
points is \p{20}{i,jk\ell}, \p{20}{i,jkm} and the period-$2$
superattractor is at \q{60}{i,jk_1}, \q{60}{i,jk_2}. Once again, at
these points $h_{11}$ is repelling off the line. Dynamical experiments
on the respective lines show that these points attract all six critical
points. Thereby, the associated basins have full measure \emph{on} the
lines. Basin portraits for these restrictions appear in
\figs~\ref{fig:h11L15} and \ref{fig:h11L30}. Since almost every point on
these these lines is in the basin of an overall saddle point, the lines
themselves behave as saddles and, thereby, are measure-zero pieces of the 
Julia set $J_{h_{11}}$.

Since the pair of $15$-lines \LL{1}{15_{ij,k\ell}} and
\MM{1}{15_{ij,k\ell}} are pointwise fixed by the involution 
$$x_i \longleftrightarrow x_j\ ,\ x_k \longleftrightarrow x_\ell,$$ 
a \G{120}-equivariant that does not smash down \LL{2}{10_{ij}}
permutes these lines as sets. 

\begin{fact}

Under $h_{11}$, \MM{1}{15_{ij,k\ell}} maps to itself. With the pair of
$30$-points \q{30}{m,ij_1}, \q{30}{m,ij_2} at $0$ and $\infty$, 
$$
h_{11}|_{\MM{1}{15_{ij,k\ell}}}: 
 z \longrightarrow \frac{z(z^2 + 6)}{6\,z^2 + 1}. 
$$

\end{fact} \noindent
This map has non-critical, attracting fixed points at \p{10}{ij_1},
\p{10}{k\ell_1}. Since the four critical points belong to the associated
basins, the dynamics on the line is reliable. Also passing through the
attracting $10$-point \p{10}{ij_1} are three $10$-lines
\LL{1}{10_{k\ell}} ($k, \ell \neq i, j$) so that, at this point, $h_{11}$
repels away from the line. Hence, this line also lies in the Julia set.

The special geometry of $h_{11}$ forces a number of points to blow up:
$$
\p{5}{i},\ \p{10}{ij_1},\ \p{10}{ij_2},\ \p{15}{i,jk},\ \p{20}{i,jk\ell},\ 
\p{30}{ij,k\ell},\ \q{30}{i,jk_-},\ \q{24}{ijk\ell}.
$$
Experimental evidence suggests that neighborhoods of these blown-up
points are filled by basins of the octahedral $20$-points. Indeed, the
\CP{2} of directions through a $10$-point \p{10}{ij_1} maps to the point
itself. Lying at the intersection of three $10$-lines \LL{1}{10_{k\ell}}
($k,\ell \neq i,j$), such a location might be called
\emph{super-repelling}. In contrast, the directions through a $30$-point
\p{30}{ij,k\ell} blow up onto the superattracting $10$-line
\LL{1}{10_{ij}} whose ``basin" is that of the $20$-points
\q{20}{ij_-}.

Since the coefficients of $h_{11}$ are real, the map also preserves
\RR---the \Sym{5}-symmetric \RP{3}---as well as the \RP{2} intersections
of \RR\ with \LL{2}{5_i} and \LL{2}{10_{ij}}. In the former case there
are four \RP{1} intersections of the \RP{2} with the $10$-lines
\LL{1}{10_{ij}} while in the latter there is a single such intersection.
The stabilizer of the respective $5$-plane or $10$-plane fixes its
resident \RP{1}s. Thus, each such \RP{1} is an ``equatorial slice" of
the associated \CP{1}. Being equivalent to the map
$$z \longrightarrow -\frac{1}{z^2}$$
on the unit circle $\{|z|=1\}$, $h_{11}$ acts chaotically when
restricted to such a slice. Hence, each \RP{1} is a chaotic
attractor on the respective \RP{2}. A basin portrait for the $5$-plane reveals 
no basins other than those
of the four $10$-lines. (See
\fig{fig:h11R5}.) The \RP{2} dynamics on the $10$-plane shows, in addition to 
the chaotic line-attractor, three additional basins at the $30$-points
\p{30}{ij,k\ell}. (See \fig{fig:h11R10}.) A $30$-point belongs to the $10$-line 
\LL{1}{10_{k\ell}}, which intersects the $10$-plane \LL{2}{10_{ij}} 
transversely. Thus, in a 
neighborhood of the $30$-point, but off the $10$-plane, there is only the 
``pipe-basin" of the $20$-points \q{20}{ij_-}. Hence, the basins on the 
$10$-plane are 2-dimensional.

\begin{conj} \label{conj:h11}

The $20$-point orbit is the attractor for $h_{11}$ and the corresponding
basins have full measure in \CP{3}.

\end{conj}

Iteration experiments on \RR\ reveal attraction only to the ten
chaotically attracting \RP{1} intersections $\RR \cap \LL{1}{10_{ij}}$.

\begin{conj}

The \Sym{5}-invariant \RP{3} is non-attracting (repelling?) and so
belongs to $h_{11}$'s Julia set.

\end{conj}

\subsection{What to look for in an attractor}

A pair of $20$-points \q{20}{ij_1}, \q{20}{ij_2} associates
canonically with an orbit of \emph{ten} lines. However, there is no such
correspondence between a pair of $20$-points and an orbit of size five;
the $20$-points do not decompose into five sets of four \Sym{4} orbits.
An association of this kind makes for a natural solution to the quintic.
What could serve the purpose better than a map whose attractor is the
$5$-point orbit?

%----------

\subsection{A special map in degree six} \label{sec:6map}

In the configuration of $10$-lines \MM{1}{10_{ijk}} each $5$-point lies
at the intersection of four lines. (See
Section~\ref{sec:config}.) Moreover, these are the only intersections of
$10$-lines. To take advantage of this structure, a map could have
superattracting pipes along the $10$-lines and basins of attraction at
the $5$-points.

The family of $6$-maps has (homogeneous) dimension six. Obtaining maps for which the
$10$-lines are critical in the ``off-line" directions uses four
parameters. For the remaining two, we get a map $f_6$ whose restriction
to a $10$-line \MM{1}{ijk} is
$$z \longrightarrow z^4$$
in coordinates where the $5$-points \p{5}{\ell}, \p{5}{m} on \MM{1}{ijk}
are $0$ and $\infty$.  In hyperplane coordinates,
\small
\begin{align*}
\phi_6 =&\ 
2\,(9\,\Phi_2\,\Phi_3 - 10\,\Phi_5)\,\phi_1 - 2\,(\Phi_2^2 -
5\,\Phi_4)\,\phi_2 + 20\,\Phi_3\,\phi_3 + 15\,\Phi_2\,\phi_4\\
=&\
\bigl[ 
2\,u_1^6 - 4\,u_1\,u_2^5 - 74\,u_1^2\,u_2^3\,u_3 - 46\,u_1^3\,u_2\,u_3^2
- 14\,u_2^2\,u_3^4 - 2\,u_1\,u_3^5 - 38\,u_1^3\,u_2^2\,u_4 \\
&\
- 44\,u_1^4\,u_3\,u_4 - 50\,u_2^3\,u_3^2\,u_4 -
122\,u_1\,u_2\,u_3^3\,u_4 - 14\,u_2^4\,u_4^2 - 152\,u_1\,u_2^2\,u_3\,
u_4^2 \\
&\ 
- 68\,u_1^2\, u_3^2\,u_4^2 - 72\,u_1^2\,u_2\,u_4^3 - 22\,u_3^3\,u_4^3 -
29\,u_2\,u_3\,u_4^4 - u_1\,u_4^5,\\[5pt]
&\ -2\,u_1^5\,u_2 + 2\,u_2^6 - 44\,u_1\,u_2^4\,u_3 - 68\,u_1^2\,u_2^2\,
u_3^2 - 22\,u_1^3\, u_3^3 - u_2\,u_3^5\\
&\ - 46\,u_1^2\,u_2^3\,u_4 - 122\,u_1^3\,u_2\,u_3\,u_4 -
72\,u_2^2\,u_3^3\,u_4 - 29\,u_1\,u_3^4\,u_4 - 14\,u_1^4\,u_4^2 \\
&\ - 38\,u_2^3\,u_3\,u_4^2- 152\,u_1\,u_2\,u_3^2\, u_4^2 -
74\,u_1\,u_2^2\, u_4^3 - 50\,u_1^2\,u_3\, u_4^3 - 14\,u_3^2\, u_4^4 -
4\,u_2\,u_4^5,\\[5pt]
&\ -14\,u_1^4\,u_2^2 - 4\,u_1^5\,u_3 - u_2^5\,u_3 -
72\,u_1\,u_2^3\,u_3^2 - 38\,u_1^2\,u_2\,u_3^3 + 2\,u_3^6 -
29\,u_1\,u_2^4\,u_4 \\
&\ - 152\,u_1^2\,u_2^2\,u_3\, u_4 - 74\,u_1^3\,u_3^2\, u_4 -
44\,u_2\,u_3^4\,u_4 - 50\,u_1^3\,u_2\,u_4^2 - 68\,u_2^2\,u_3^2\, u_4^2\\
&\ - 46\,u_1\,u_3^3\, u_4^2 - 22\,u_2^3\, u_4^3- 122\,u_1\,u_2\,u_3\,
u_4^3 - 14\,u_1^2\, u_4^4 - 2\,u_3\,u_4^5,\\[5pt]
&\ -22\,u_1^3\,u_2^3 - 29\,u_1^4\,u_2\,u_3 - 14\,u_2^4\,u_3^2 -
50\,u_1\,u_2^2\,u_3^3 - 14\,u_1^2\,u_3^4 - u_1^5\,u_4 - 2\,u_2^5\,u_4\\
&\ - 122\,u_1\,u_2^3\,u_3\,u_4 - 152\,u_1^2\,u_2\,u_3^2\, u_4 -
4\,u_3^5\,u_4 - 68\,u_1^2\,u_2^2\, u_4^2 - 72\,u_1^3\,u_3\, u_4^2 \\
&\ - 74\,u_2\,u_3^3\, u_4^2 - 46\,u_2^2\,u_3\, u_4^3 - 38\,u_1\,u_3^2\,
u_4^3 - 44\,u_1\,u_2\, u_4^4 + 2\,u_4^6 
\bigr]. \end{align*} \normalsize

By construction, $f_6$ self-maps each \Sym{3}-symmetric $10$-plane
\LL{2}{10_{ij}}. The $10$-point \p{10}{ij_2} and $5$-points \p{5}{k} ($k
\neq i,j$) form \Sym{3} orbits on \LL{2}{10_{ij}} of sizes one and
three. Furthermore, $f_6$ preserves \RR---the \Sym{5}-symmetric \RP{3}.
We can get a picture of
the map's \emph{restricted dynamics} by plotting basins of attraction on
the \RP{2} intersection
$$\LL{2}{10_{ij}} \cap \RR.$$
(See \figs~\ref{fig:f6R10},
through \ref{fig:f6R10mag2} in \app{app:basins}.)
The plot shows attraction to the $5$-points and the $10$-point. However,
the $10$-point lies on the ``equator" of an \MM{1}{10_{k\ell m}}
($k,\ell,m \neq i,j$) where $f_6$ \emph{repels} in the off-plane
direction. Thus, the $2$-dimensional basin of a $10$-point is a measure-zero
part of $J_{f_6}$.  No other attracting sets appear. Moreover, regions of
positive measure that do not belong to one of these four ``restricted
basins" are not evident. The plot is consistent with the claim that the
only fully $3$-dimensional basins are those of the $5$-points.

A $15$-line \LL{1}{15_{ij,k\ell}} contains one $5$-point \p{5}{m}, one
$15$-point \p{15}{m,ij} ($m \neq k,\ell$), and two $10$-points
\p{10}{ij_2}, \p{10}{k\ell_2}. In coordinates where the $5$-point is
$0$, the $15$-point is $\infty$, and the $10$-points are $\pm 1$ the map
restricts to
$$z \longrightarrow \frac{48\,z^5}{-3 - z^2 + 35\,z^4 + 17\,z^6}\,.$$
The critical points of the restricted map are 
$$0,\ \pm 1,\ \pm \sqrt{\frac{9 \pm 4\,\sqrt{21}}{17}}\,.$$
with $0,\pm 1$ fixed. Experiment reveals that the four non-fixed
critical points belong to the basins of the three superattracting
points. Hence, these basins have full measure on the
$15$-line. (\figs~\ref{fig:f6L15} and \ref{fig:f6L15mag} display
portraits.) As a member of three $15$-lines \LL{1}{15_{ij,k\ell}} a
$10$-point \p{10}{ij_2} superattracts in these directions. However,
these three lines lie in the $10$-plane \LL{2}{10_{ij}} so that, as seen
above, $f_6$ is completely superattracting \emph{in} the plane at
\p{10}{ij_2}.

Another distinction for $f_6$ is its action on a $15$-line
\MM{1}{15_{ij,k\ell}} which, by equivariance, must map either to itself
or \LL{1}{15_{ij,k\ell}}.

\begin{fact}

Under $f_6$, \MM{1}{15_{ij,k\ell}} maps to \LL{1}{15_{ij,k\ell}}.
Effectively, this creates a second orbit of superattracting pipes to the
$5$-points.

\end{fact} \noindent
This is what led me to $6$-maps, each of which send the $10$-point
\p{10}{ij_1} to the $10$-point \p{10}{ij_2}.

Finally, noting that $\phi_6$ has real coefficients, it must preserve
the \RP{3} whose points have real $u$ coordinates. This is
\emph{not} the \Sym{5}-symmetric \RR. Rather it seems to be associated
with the \Sym{4} stabilizer of \p{5}{1} which is $[1,1,1,1]$ in
the $u$ space. This \RP{3} intersects the $10$-planes \LL{2}{10_{25}} and
\LL{2}{10_{34}} in an \RP{2} with $\Z{2} \times \Z{2}$ symmetry. In
addition to \p{5}{1} this \RP{2} contains the $10$-points \p{10}{25_1},
\p{10}{25_2}, \p{10}{34_2} as well as the \RP{1} through \p{10}{25_1}
and \p{10}{25_2}. Since this line is an equatorial slice through
\LL{1}{10_{125}}, $f_6$ attracts chaotically along the line. (See
\fig{fig:f6R5L10} for a basin portrait.)

Graphical and experimental evidence supports the claim of reliability
for $f_6$.

\begin{conj} \label{conj:f6}

The attractor for $f_6$ is the $5$-point orbit the basins of which fill
up \CP{3} in measure.

\end{conj}

%-------

\section{Solving the quintic} \label{sec:quintic}

To compute a root of a polynomial, one must overcome the symmetry
present. For a general equation of degree $n$ the obstacle is \Sym{n}.
Klein described a means to this end: given values for an ``independent"
set of \Sym{n}-invariant homogeneous polynomials
$$a_1=G_1(x),\ \dots,\ a_m=G_m(x),$$
find the \Sym{n} orbits of solutions $x$ to these
equations. \cite[pp. 69ff]{Klein} This task of inverting the
$G_k$ is the \emph{form problem on \Sym{n}}. It also has a rational
manifestation: for $m-1$ given values, invert $m-1$ invariant rational
functions of degree zero.

An \Sym{n} equivariant with reliable dynamics breaks the obstructing
symmetry. In effect, this provides a mechanism for solving the form
problem and, hence, the $n$th degree equation.  What follows is one way to use
\G{120}-symmetry in multiple settings to assemble a procedure that solves
almost any quintic.

%----------

\subsection{Parameters}

The \G{120} rational form problem is to solve 
\begin{equation} \label{eq:K}
K_1 = \frac{\Phi_4(u)}{\Phi_2(u)^2} \quad
K_2 = \frac{\Phi_3(u)^2}{\Phi_2(u)^3} \quad 
K_3 = \frac{\Phi_5(u)}{\Phi_2(u)\,\Phi_3(u)}. 
\end{equation} 
As functions, the $K_i$ define the \G{120} quotient map
$$
[K_1,K_2,K_3,1] = 
[\Phi_2\,\Phi_3\,\Phi_4,\Phi_3^3,\Phi_2^2\,\Phi_5,\Phi_2^3\,\Phi_3]
$$
on $\CP{3} \setminus \{\Phi_2=\Phi_3=0\}$. The generic fiber over points
in \CP{3} is a \G{120} orbit given by
$$
\{\Phi_4 = 
K_1\,\Phi_2^2\} \cap \{\Phi_3^2 = K_2\,\Phi_2^3\} \cap 
\{\Phi_5 = K_3\,\Phi_2\,\Phi_3\}.
$$
Exceptional locations are $[0,1,0,0]$ and $[0,0,1,0]$ where the
respective fibers include the quadric and cubic surfaces $\{\Phi_2=0\}$
and $\{\Phi_3=0\}$.

Between quintic equations and \G{120} actions the parameters $K_i$ forge
a link. The connection consists in $K$-parametrizations of each regime.  From a
parametrized family of \G{120} actions, we can extract parametrized families of
\Sym{5} invariants and equivariant $6$-maps.  In this way, a choice of parameter
$K$ produces a quintic $R_K$ as well as a system of invariants $\Phi_{2_K}(w),
\dots, \Phi_{5_K}(w)$, and a $6$-map $\phi_K(w)$---a conjugate of $\phi_6(u)$---
on a parametrized $w$-space. 

%----------

\subsection{A family of \Sym{5} quintics}

Let \G{v} be a version of \G{120} that acts on a $v$-coordinatized
$\CP{3}_v$. This will be a parameter space---the coordinate $v$ merely
stands in for $u$. The linear polynomials
$$X_{k}(x) = - 4\,x_k + \sum_{i \neq k} x_i $$
form an orbit of size five.  In hyperplane coordinates, these are
\begin{align*}
L_1(v) &= -\sqrt{5}\,(v_1 + v_2 + v_3 + v_4)\\
L_2(v) &= -\sqrt{5}\,\w{5}\,(\w{5}^3\,v_1 + \w{5}^2\,v_2 + \w{5}\,v_3 +
v_4)\\
L_3(v) &= -\sqrt{5}\,\w{5}\,(\w{5}^2\,v_1 + v_2 + \w{5}^3\,v_3 +
\w{5}\,v_4) \\
L_4(v) &= -\sqrt{5}\,\w{5}\,(\w{5}\,v_1 + \w{5}^3\,v_2 + v_3 +
\w{5}^2\,v_4)\\
L_5(v) &= -\sqrt{5}\,\w{5}\,(v_1 + \w{5}\,v_2 + \w{5}^2\,v_3 +
\w{5}^3\,v_4).
\end{align*}
The rational functions
$$S_k(v) = \frac{\Phi_2(v)\,L_{k}(v)}{\Phi_3(v)} $$
also give a 5-orbit.  Taking the $S_k$ as roots of a polynomial
$$
R_v(s) = \prod_{k=1}^5 \bigl(s - S_k(v)\bigr) = 
\sum_{k=0}^5 C_k(v)\,s^{5-k}
$$
yields a family of quintics whose members generically have \Sym{5}
symmetry. Since \G{v} permutes the $S_k(v)$, each coefficient $C_k(v)$
is \G{v}-invariant and hence, expressible in terms of the basic forms
$\Phi_k(v)$ and, ultimately, in terms of the $K_i$.  Of course, $C_0(v)=1$.  
Since there is no degree-$1$ \G{120} invariant, $C_1(v)=0$.  Direct calculation
determines the remaining coefficients:
\small \begin{align*}
C_2(v) =&\ -\frac{125\,\Phi_2(v)^3}{2\,\Phi_3(v)^2} = 
-\frac{125}{2\,K_2(v)}\\
C_3(v) =&\ \frac{625\,\sqrt{5}\,\Phi_2(v)^3}{3\,\Phi_3(v)^2} =
           \frac{625\,\sqrt{5}}{3\,K_2(v)}\\
C_4(v) =&\ 
\frac{15625\,\bigl(\Phi_2(v)^6 - 2\,\Phi_2(v)^4 \Phi_4(v)\bigr)}{8\,\Phi_3(v)^4} 
= \frac{15625\,\bigl(1 - 2\,K_1(v)\bigr)}{8\,K_2(v)^4}\\
C_5(v) =&\ \frac{-15625\,\bigl( 5\,\sqrt{5}\,\Phi_2(v)^6\,\Phi_3(v) - 
          6\,\sqrt{5}\,\Phi_2(v)^5\,\Phi_5(v)\bigr)}{6\,\Phi_3(v)^5} =
           \frac{-15625\,\sqrt{5}\,( 6\,K_3(v) - 5)}{6\,K_2(v)^2}.
\end{align*} \normalsize
Members of the 3-parameter family of quintic \G{120} \emph{resolvents}
$$
R_K(s) = s^5  - 
  \frac{125}{2\,K_2}\,s^3 + 
  \frac{625\,\sqrt{5}}{3\,K_2}\,s^2 - 
  \frac{15625\,\left( -1 + 2\,K_1 \right)}{8\,K_2^2}\,s+ 
  \frac{15625\,\sqrt{5}\,\left( -5 + 6\,K_3 \right)}{6\,K_2^2}
$$
are particularly well-suited for an iterative solution that employs
$\phi_6$. For selected values of the $K_i$, a solution to the resulting
form problem yields a root of $R_K$. Use of \G{120} symmetry will provide a
means of finding such a solution without explicitly inverting the $K_i$
equations (\ref{eq:K}).

%----------

\subsection{Reduction of the general quintic to a \G{120} resolvent}

By means of a well-known linear \emph{Tschirnhaus} transformation
the general quintic
becomes the \emph{standard} $4$-parameter \emph{resolvent}
$$q(y) = y^5 + b_2\,y^3 + b_3\,y^2 + b_4\,y + b_5.$$

Application of another linear \emph{Tschirnhaus} transformation
$$s \longrightarrow \frac{y}{\lambda}$$
converts the $3$-parameter family $R_K(s)$ into a \G{120}
resolvent
\begin{align*}
\Sigma_{K,\lambda}(y) =&\ \lambda^5\,R_K\Bigl(\frac{y}{\lambda}\Bigr) 
=\ y^5 + \lambda^2\,C_2\,y^3 + \lambda^3\,C_3\,y^2
+ \lambda^4\,C_4\,y + \lambda^5\,C_5
\end{align*}
in the four parameters $K_1,K_2,K_3$, and the \emph{auxiliary}
$\lambda$.

The functions 
$$b_k = \lambda^k\,C_k$$
relate the coefficients of $q$ and $\Sigma_{K,\lambda}$. These invert to
$$
K_1 = \frac{b_2^2 - 2\,b_4}{2\,b_2^2}\qquad
K_2 =\ \frac{-9\,b_3^2}{8\,b_2^3}\qquad
K_3 = \frac{5\,\bigl(b_2\,b_3 - b_5 \bigr)}{6\,b_2\,b_3}\qquad
\lambda = \frac{-3\,b_3}{10\,\sqrt{5}\,b_2}.
$$
Thus, almost any quintic descends to a member of $R_K$. The reduction
fails when
$$
-2\,a_1^2 + 5\,a_2 = 5\,b_2 = 0 \qquad \text{or} \qquad
4\,a_1^3 - 15\,a_1\,a_2 + 25\,a_3 = 25\,b_3 = 0.
$$
A solution to the special resolvent $R_K$ then ascends to a solution to
the general quintic.

%----------

\subsection{A family of \Sym{5} actions}

With the basic \G{v}-maps, construct the \emph{parametrized change of
coordinates}
$$
u = \tau_v w = \sum_{i=1}^4 \bigl(\Phi_{6-i}(v)\,\phi_i(v) \bigr) w_i.
$$
A matrix form results from taking the $\phi_k(v)$ as column vectors:
\small $$
\begin{pmatrix} u_1\\u_2\\u_3\\u_4 \end{pmatrix} = \begin{pmatrix}
&&&&&&\\
\Phi_5(v)\,\phi_1(v)&\vdots&
\Phi_4(v)\,\phi_2(v)&\vdots&
\Phi_3(v)\,\phi_3(v)&\vdots&
\Phi_2(v)\,\phi_4(v)\\
&&&&&&
\end{pmatrix}
\begin{pmatrix} w_1\\w_2\\w_3\\w_4 \end{pmatrix}.
$$ \normalsize
For a choice of parameter $v$,
$$\tau_v: \CP{3}_w \longrightarrow \CP{3}_u$$
is linear in $w$ and gives rise to a parametrized family of \G{120} groups
$$\G{w}^v = \tau_v^{-1} \G{u} \tau_v.$$
The setup here is as follows.
\begin{itemize}

\item \G{u} is a version of \G{120} that acts on a \emph{reference space}
$\CP{3}_u$.

\item \G{v} is a version of \G{120} that acts on a \emph{parameter space}
$\CP{3}_v$.

\item \G{u} and \G{v} have identical expressions in their respective
coordinates.

\item $\G{w}^v$ is a version of \G{120} that acts on a
\emph{parametrized space} $\CP{3}_w$.

\item The iteration that solves quintics in $R_K$ will take place in $\CP{3}_w$.

\end{itemize}

Each $\G{w}^v$ has its system of invariants and equivariants. From this
point of view, we can see, in the resolvents $R_v$ and $\G{w}^v$
equivariants, a connection between quintics and dynamical systems.
Furthermore, each $\G{w}^v$ invariant and equivariant is expressible in
the $K_i$.

The first thing to notice is that, by construction, $\tau_v w$ possesses
an equivariance property:
$$
\tau_{A v} w = A\,\tau_v w \quad \text{for}\ A \in \G{v}, \G{u}.
$$
The determinant of $\tau_v$ will enter into upcoming calculations and so,
demands some attention. Since
\begin{equation} \label{eq:tauDet}
|\tau_{A v}| = |A|\,|\tau_v|,
\end{equation}
$|\tau_v|$ is invariant under the \Alt{5} subgroup \G{60} of \G{v} but
only relatively invariant under the full \Sym{5} group \G{120}. The even
transformations have determinant $1$ while the odd elements have
determinant $-1$. Furthermore, 
\begin{align*} 
|\tau_v| =&\
\Phi_2(v)\,\Phi_3(v)\,\Phi_4(v)\,\Phi_5(v)\, 
\begin{vmatrix} 
\phi_1(v)\ |\ \phi_2(v)\ |\ \phi_3(v)\ |\ \phi_4(v)
\end{vmatrix}\\ 
=&\ \Phi_2(v)\,\Phi_3(v)\,\Phi_4(v)\,\Phi_5(v)\,\Psi_{10}(v) 
\end{align*}
where $\Psi_{10}$ is a scalar multiple of the product of the ten linear
forms associated with the ten planes of reflection that generate \G{v}.
Reflection group theory tells us that this is the only form in degree ten that 
is invariant under \G{60} but not \G{120}. From (\ref{eq:tauDet}), the 
degree-$48$ square of $|\tau_v|$ is \G{120}-invariant. Let
$$|\tau_v|^2 = \Phi_2^{24}(v)\,t_K$$
determine its $K$-expression. The explicit form of $t_K$ appears in
\app{app:forms}.

%----------

\subsection{A family of \Sym{5} invariants}

The equivariance in $v$ of $\tau_v w$ implies that $\Phi_2(\tau_v w)$ is
\G{v}-invariant. Thus, each $w$ coefficient of $\Phi_2(\tau_v w)$ inherits
the same invariance. Since
$$
\deg_v \Phi_2(\tau_v w) = \deg_u \Phi_2(u) \cdot \deg_v \tau_v w = 2 \cdot 6 = 12,
$$ 
the rational function
$$\frac{\Phi_2(u)}{\Phi_2(v)^6} = \frac{\Phi_2(\tau_v w)}{\Phi_2(v)^6}$$ 
is degree zero in $v$ and thereby, expressible in the $K_i$. Let 
\begin{equation} \label{eq:Phi2}
\Phi_2(v)^6\,\Phi_{2_K}(w) = \Phi_2(u)
\end{equation}
define the basic degree-2 $\G{w}^v$ invariant $\Phi_{2_K}(w)$. Solving a
system of linear equations whose dimension is that of the degree-12
\G{v} invariants yields an explicit expression in
the $K_i$ for each $w$-coefficient of $\Phi_2(\tau_v w)$. Similar
considerations apply in degree three where
\begin{equation} \label{eq:Phi3}
\Phi_2(v)^9\,\Phi_{3_K}(w) = \Phi_3(u).
\end{equation}
The results appear in \app{app:forms}.

By Fact~\ref{fact:diffInv}, the degree-$4$ and degree-$5$
invariants derive from those in degrees two and three. First of all, the
chain rule determines transformation formulas for the hessian and
bordered hessian.

\begin{prop}

For $y=A x$,
\begin{align*}
H_x(F(y)) =&\ A^T H_y(F(y)) A\\
B_x(F(y),G(y)) =&\ 
\begin{pmatrix}A^T&0\\0&1 \end{pmatrix} 
B_y(F(y),G(y)) 
\begin{pmatrix}A&0\\0&1 \end{pmatrix}
\end{align*}
where the subscript indicates the variable of differentiation.  Thus,
\begin{align*}
\bigl| H_x(F(y)) \bigr| =&\ |A|^2 \bigl| H_y(F(y)) \bigr|\\
\bigl| B_x(F(y),G(y)) \bigr| =&\ |A|^2 \bigl|B_y(F(y),G(y)) \bigr|.
\end{align*}

\end{prop} \noindent
Applied to the parametrized change of variable $w = \tau_v^{-1} u$,
\begin{align*} 
G_4(u) =&\ \bigl|H_u \bigl(\Phi_3(u)\bigr) \bigr|\\
=&\ \Bigl|H_u \bigl(\Phi_2(v)^9\,\Phi_{3_K}(w) \bigr) \Bigr|\\
=&\ \frac{\bigl(\Phi_2(v)^9 \bigr)^4}{|\tau_v|^2}\,
\Bigl|H_w \bigl(\Phi_{3_K}(w) \bigr) \Bigr|\\
=&\ \frac{\Phi_2(v)^{12}}{t_K}\,
\Bigl|H_w \bigl(\Phi_{3_K}(w) \bigr)\Bigr|\\
=&\ \Phi_2(v)^{12}\,G_{4_K}(w)
\end{align*}
and
\begin{align*} 
G_5(u) =&\ \bigl|B_u \bigl(\Phi_3(u),\Phi_2(u)\bigr) \bigr|\\
=&\ \Bigl|B_u 
 \bigl(\Phi_2(v)^9\,\Phi_{3_K}(w),\Phi_2(v)^6\,\Phi_{2_K}(w) \bigr) 
\Bigr|\\
=&\ \frac{\bigl(\Phi_2(v)^9 \bigr)^3 \bigl(\Phi_2(v)^6 \bigr)^2}
    {|\tau_v|^2}\,
   \Bigl|B_w \bigl(\Phi_{3_K}(w),\Phi_{2_K}(w) \bigr) \Bigr|\\
=&\ \frac{\Phi_2(v)^{15} }{t_K}\,
\Bigl|B_w \bigl(\Phi_{3_K}(w),\Phi_{2_K}(w) \bigr) \Bigr|\\
=&\ \Phi_2(v)^{15}\,G_{5_K}(w).
\end{align*}
Employed here are the obvious definitions
$$
G_{4_K}(w) = \frac{\Bigl|H_w \bigl(\Phi_{3_K}(w) \bigr) \Bigr|}{t_K}
\quad
G_{5_K}(w) = \frac{
\Bigl|B_w \bigl(\Phi_{3_K}(w),\Phi_{2_K}(w) \bigr) \Bigr|}{t_K}.
$$
With natural definitions for $\Phi_{4_K}(w)$ and $\Phi_{5_K}(w)$,
\begin{align*}
\Phi_4(u) 
=&\ \frac{1}{324} \bigl(\Phi_2(u)^2 - 5\,G_4(u) \bigr)\\ 
=&\ \frac{1}{324} \bigl(
 \Phi_2(v)^{12}\,\Phi_{2_K}(w)^2 - 5\,\Phi_2(v)^{12}\,G_{4_K}(w)\bigr)\\
=&\ \Phi_2(v)^{12}\,\Phi_{4_K}(w)
\end{align*}
and
\begin{align*}
\Phi_5(u) 
=&\ \frac{1}{864} \bigl(720\,\Phi_2(u)\,\Phi_3(u) + G_5(u) \bigr)\\ 
=&\ \frac{1}{864} \bigl( 720\,
\Phi_2(v)^{15}\,\Phi_{2_K}(w)\,\Phi_{3_K}(w) 
+ \Phi_2(v)^{15}\,G_{5_K}(w) \bigr)\\
=&\ \Phi_2(v)^{15}\,\Phi_{5_K}(w).
\end{align*}

%----------

\subsection{A family of \Sym{5} equivariant $6$-maps}

Emerging from each $\G{w}^v$ action is a version
$\tau_v^{-1}\,\phi_6(\tau_v w)$ of $\phi_6(u)$. Being \G{v}-invariant,
these maps also admit parametrization by $K$. Thereby, each quintic
$R_K$ enters into association with a dynamical system $\phi_K$ on
$\CP{3}_w$.

The reversed identity $R$ and gradient $\nabla^r=R\,\nabla$ appeared in the
context of a change from five $x$ coordinates to four $u$ coordinates. In the
present setting, a \emph{reversed transpose} arises.
\begin{defn}

The \emph{repose} $A^r$ of an $n \times n$ matrix $A$ is its reflection
through the \emph{reversed diagonal}---the entries whose subscripts sum
to $n+1$. Alternatively, 
$$A^r = R\, A^T\,R.$$

\end{defn}

\begin{prop}

For a change of coordinates $u=A\,w$ and a polynomial $\Phi(u) =
\tilde{\Phi}(w)$, the reversed gradient map transforms by
$$\nabla_u^r \Phi(u) = A^r\,\nabla_w^r \tilde{\Phi}(w).$$
\end{prop}

\begin{proof}

Noting that $R^2=I$,
\begin{align*}
\nabla_u^r \Phi(u) =&\ R\,\nabla_u \Phi(u) \\
=&\ R\,A^T\,\nabla_w \tilde{\Phi}(w) \\
=&\ R\,A^T\,R\,R\,\nabla_w \tilde{\Phi}(w) \\
=&\ A^r\,\nabla_w^r \tilde{\Phi}(w).
\end{align*}

\end{proof}
For the generating \G{120} maps,
\begin{align*}
\phi_\ell(u) =&\ \nabla_u^r \Phi_2(v)^{3(\ell + 1)}\,\Phi_{{\ell+1}_K}(w) \\
=&\ \Phi_{2}(v)^{3 (\ell + 1)}\,
 (\tau_v^{-1})^r\,\nabla_w^r \Phi_{{\ell+1}_K}(w) \\
=&\ \Phi_{2}(v)^{3 (\ell + 1)}\,\tau_v\,\tau_v^{-1}\,(\tau_v^{-1})^r\,
   \nabla_w^r \Phi_{{\ell+1}_K}(w) \\
=&\ \tau_v\,\Phi_{2}(v)^{3 (\ell + 1)}\,(\tau_v^r\,\tau_v)^{-1}\,
   \nabla_w^r \Phi_{{\ell+1}_K}(w).
\end{align*}
Thus,
$$
\tau_v^{-1}\,\phi_\ell(\tau_v w) = 
\Phi_{2}(v)^{3 (\ell + 1)}\,
(\tau_v^r\,\tau_v)^{-1}\,\nabla_w^r \Phi_{{\ell+1}_K}(w).
$$
Using the description on the left-hand side, a straightforward
calculation reveals this map to be invariant in $v$ so that the matrix
$\tau_v^r\,\tau_v$ has entries that are degree-$12$ \G{v} invariants.
Hence, the matrix product has a $K$-expression:
$$
\tau_v^r\,\tau_v = \Phi_2(v)^6\,T_K 
\quad \text{or} \quad
(\tau_v^r\,\tau_v)^{-1} = \frac{T_K^{-1}}{\Phi_2(v)^6}.
$$
(See \app{app:forms} for the explicit form.)
Using this to express the transformation of basic equivariants yields
\begin{align*}
\phi_\ell(u) =&\ \Phi_{2}(v)^{3 (\ell + 1)}\,\tau_v\,
  \frac{T_K^{-1}}{\Phi_2(v)^6}\,
   \nabla_w^r \Phi_{{\ell+1}_K}(w) \\
=&\ \Phi_{2}(v)^{3 (\ell-1)}\,\tau_v\,\phi_{\ell_K}(w)
\end{align*}
where 
$$\phi_{\ell_K}(w) = T_K^{-1}\,\nabla_w^r \Phi_{{\ell+1}_K}(w).$$

Finally, we can identify a $K$-parametrized $6$-map $\phi_K(w)$ that is
conjugate to $\phi_6(u)$. The map's expression in basic terms appears after 
substitution into the
formula found in Section~\ref{sec:6map}.  (See
\app{app:forms}.)

%----------

\subsection{Root selection}

Being conjugate to $\phi_6(u)$ each $\phi_K(w)$ shares the former's
conjectured reliable dynamics. Accordingly, the attractor for each
choice of $K_i$ is the $5$-point orbit in the corresponding $\CP{3}_w$
so that for almost every $w_0 \in \CP{3}_w$,
$$
\phi_K^n(w_0) \longrightarrow \tau_v^{-1} \p{5}{\ell} \qquad
\text{for some $5$-point}\ \p{5}{\ell} \in \CP{3}_u.
$$
To solve the resolvent $R_K$, the output of the iteration must link with
the roots of $R_K$. With this, we see that solving $R_K$ amounts to
inverting $\tau_v$---the form problem in yet another guise. With the
assistance of a \G{120} tool, this is effectively what the dynamics of
$\phi_K$ accomplishes. (This clever device is due to McMullen.)

To manufacture the root-selecting tool, we begin with an orbit of quadratic 
\Sym{4}-invariants
$$X^2_{k}(x) = - 4\,x_k^2 + \sum_{i \neq k} x_i^2. $$
These form a \G{120} orbit of size five.  Their hyperplane expressions are
\small \begin{align*}
Q_1(u) =&\ 
2\,\bigl( 2\,u_1^2 + 4\,u_1\,u_2 + 2\,u_2^2 + 4\,u_1\,u_3 + 3\,u_2\,u_3 +
2\,u_3^2 + 3\,u_1\,u_4 \\
&\ + 4\,u_2\,u_4 + 4\,u_3\,u_4 + 2\,u_4^2 \bigr)\\
Q_2(u) =&\
2\,\bigl( 2\,\w{5}^3\,u_1^2 + 4\,\w{5}^2\,u_1\,u_2 + 2\,\w{5}\,u_2^2 +
4\,\w{5}\,u_1\,u_3 + 3\,u_2\,u_3 - 2\,u_3^2 - 2\,\w{5}\,u_3^2 \\
&\ - 2\,\w{5}^2\,u_3^2 - 2\,\w{5}^3\,u_3^2 + 3\,u_1\,u_4 - 4\,u_2\,u_4 -
4\,\w{5}\,u_2\,u_4 - 4\,\w{5}^2\,u_2\,u_4 \\
&\ - 4\,\w{5}^3\,u_2\,u_4 + 4\,\w{5}^3\,u_3\,u_4 
 + 2\,\w{5}^2\,u_4^2 \bigr)\\
Q_3(u) =&\ 
-2\,\bigl( -2\,\w{5}\,u_1^2 + 4\,u_1\,u_2 + 4\,\w{5}\,u_1\,u_2 +
4\,\w{5}^2\,u_1\,u_2 + 4\,\w{5}^3\,u_1\,u_2 - 2\,\w{5}^2\,u_2^2 \\
&\ - 4\,\w{5}^2\,u_1\,u_3 - 3\,u_2\,u_3 - 2\,\w{5}^3\,u_3^2 -
3\,u_1\,u_4 - 4\,\w{5}^3\,u_2\,u_4 - 4\,\w{5}\,u_3\,u_4 \\
&\ + 2\,u_4^2 + 2\,\w{5}\,u_4^2 + 2\,\w{5}^2\,u_4^2 + 2\,\w{5}^3\,u_4^2
\bigr)\\
Q_4(u) =&\ 
-2\,\bigl( 2\,u_1^2 + 2\,\w{5}\,u_1^2 + 2\,\w{5}^2\,u_1^2 + 2\,\w{5}^3\,u_1^2
- 4\,\w{5}\,u_1\,u_2 - 2\,\w{5}^3\,u_2^2 \\
&\ - 4\,\w{5}^3\,u_1\,u_3 - 3\,u_2\,u_3 - 2\,\w{5}^2\,u_3^2 - 3\,u_1\,u_4 -
4\,\w{5}^2\,u_2\,u_4 + 4\,u_3\,u_4 \\
&\ + 4\,\w{5}\,u_3\,u_4 + 4\,\w{5}^2\,u_3\,u_4 + 4\,\w{5}^3\,u_3\,u_4 -
2\,\w{5}\,u_4^2 \bigr)\\
Q_5(u) =&\ 
2\,\bigl( 2\,\w{5}^2\,u_1^2 + 4\,\w{5}^3\,u_1\,u_2 - 2\,u_2^2 - 2\,\w{5}\,u_2^2
- 2\,\w{5}^2\,u_2^2 - 2\,\w{5}^3\,u_2^2 - 4\,u_1\,u_3 \\
&\ - 4\,\w{5}\,u_1\,u_3 - 4\,\w{5}^2\,u_1\,u_3 - 4\,\w{5}^3\,u_1\,u_3 +
3\,u_2\,u_3 + 2\,\w{5}\,u_3^2 + 3\,u_1\,u_4 \\
&\ + 4\,\w{5}\,u_2\,u_4 + 4\,\w{5}^2\,u_3\,u_4 + 2\,\w{5}^3\,u_4^2 \bigr).
\end{align*} \normalsize
Furthermore, each of the five forms 
$$
G_k(u) = -\frac{3}{25}\,L_k(u)^2 + Q_k(u) \quad k = 1, \dots, 5
$$ 
vanish at the $5$-points \p{5}{\ell} with $\ell \neq k$ but not at \p{5}{k}.

Now, to draw the roots of the quintics $R_K(s)$ into the game, consider the 
rational function
$$J_v(w) = 
\alpha\,\sum_{k=1}^5 
 \frac{G_k(\tau_v w)}{\Phi_2(\tau_v w)}\,\frac{\Phi_2(v)\,L_k(v)}{\Phi_3(v)}
= \alpha\,\sum_{k=1}^5 \frac{G_k(\tau_v w)}{\Phi_2(\tau_v w)}\,S_k(v)
$$
where $\alpha$ is a constant to be determined. Since the $v$-degree of
the numerator and denominator is $15=2\cdot 6 + 3$ while the $w$-degree
is $2$, the function is rationally degree zero in both variables. At a
$5$-point $\tau_v^{-1} \p{5}{\ell}$ in $\CP{3}_w$ four of the five terms
in $J_v$ vanish; this leaves
$$\alpha\,\frac{G_\ell(\p{5}{\ell})}{\Phi_2(\p{5}{\ell})}\,S_\ell(v).$$
Setting
$$
\alpha = \frac{\Phi_2(\p{5}{1})}{G_1(\p{5}{1})} =
\dots = \frac{\Phi_2(\p{5}{5})}{G_5(\p{5}{5})} = \frac{1}{15}\ 
$$
``selects" the root $S_\ell(v)$ of $R_K(s)$. Since the iterative
``output" of $\phi_K(w)$ is a single $5$-point in $\CP{3}_w$, the
dynamics produces one root.

The root-selector $J_v(w)$ has invariance properties that allow it to exhibit a 
useful form. Let 
$$\Gamma_v(w) =\sum_{k=1}^5 G_k(\tau_v w)\,L_k(v).$$ 
Since \G{v} permutes its terms,
$\Gamma_v$ is invariant under the action and hence,
expressible in $K$:
$$\Gamma_v(w) = \Phi_2(v)^5\,\Phi_3(v)\,\Gamma_K(w).$$
(The explicit form of $\Gamma_K$ appears in
\app{app:forms}.)  Finally, application of (\ref{eq:Phi2}) yields
\begin{align*}
J_v(w) 
=&\ \frac{\Phi_2(v)\,\Gamma_v(w)}{15\,\Phi_3(v)\,\Phi_2(\tau_v w)}\\
J_K(w)=&\ \frac{\Gamma_K(w)}{15\,\Phi_{2_K}(w)}.
\end{align*}

%----------

\subsection{The procedure summarized}

\begin{enumerate}

\item Select a general $5$-parameter quintic $p(x)$.

\item Tschirnhaus transform $p(x)$ into a member $R_K(s)$ of the
$3$-parameter family of \G{120} quintics---this determines values for
$K_1, K_2, K_3$ as well as the auxiliary parameter $\lambda$.

\item For the selected $K$ values compute the invariants
$\Phi_{i_K}(w)\ (i = 2,3,4,5)$, the $6$-map $\phi_K(w)$, the form
$\Gamma_K(w)$, and the root-selector $J_K(w)$. (In fact, a
rather lengthy once-and-for-all expression for $\phi_K(w)$ is easy to
compute.  \cite{CrassWeb} Such a formula renders calculations of 
$\Phi_{2_K}$, $\Phi_{3_K}$, and $\Phi_{4_K}$ superfluous .)

\item From an arbitrary initial point $w_0$ iterate $\phi_K$ until
convergence:
$$\phi_K^n(w_0) \longrightarrow w_\infty.$$
Conjecturally, the output $w_\infty$ is a $5$-point in $\CP{3}_w$.

\item Compute a root $S = J_K(w_\infty)$ of $R_K$.

\item Transform $S$ into a root of $p(x)$.

\end{enumerate}
(At \cite{CrassWeb}, there are \emph{Mathematica}
data files and a notebook that implement the iterative solution to the
quintic.)

%--------
\clearpage
\appendix

%-------

\section{Special orbit data} \label{app:tables}

For ease of reference, the following tables provide descriptions of the
special \G{120} orbits that bear upon the quintic-solving algorithm.

\small
%%%%%%%%%%%%%
\begin{table}[hb]
 
\caption{Special points on
$\left\{[x_1,x_2,x_3,x_4,x_5]\ |\ x_k \in \R{} \right\} \simeq \RP{3}$}

\label{tab:RP3}

$$\begin{array}{c|c|c|c}

\text{Size}&\text{Representative}&
\text{Descriptor}&\text{Stabilizer}\\
\hline

&&&\\[-10pt]

5&[-4,1,1,1,1]&\p{5}{1}&\Sym{4}\\[5pt]

10&[0,0,0,1,-1]&\p{10}{45_1}&\Sym{3} \times \Z{2}\\[\tabGap]

10&[2,2,2,-3,-3]&\p{10}{45_2}&\Sym{3} \times \Z{2}\\[\tabGap]

15&[0,1,1,-1,-1]&\p{15}{1,23}=\p{15}{1,45}&\D{4}\\[\tabGap]

20&[0,-3,1,1,1]\ [-3,0,1,1,1]& \p{20}{1,345}\ \p{20}{2,345}&\Sym{3}\\[\tabGap]

30&[0,0,1,1,-2]&\p{30}{12,34}&\Z{2} \times \Z{2}\\[\tabGap]

\end{array}$$

\end{table}
%%%%%%%%%%%%%
\begin{table}

\caption{Special points on $\QQ{}{}=\{\sum_{k=1}^5 x_k^2 = 0\}$}

\label{tab:Q}

$$\begin{array}{c|c|c|c|l}

\text{Size}&\text{Representative}&
\text{Descriptor}&\text{Stabilizer}&\text{Remarks}\\
\hline

&&&\\[-10pt]

20&[0,0,1,\w{3},\w{3}^2]&\q{20}{12_1}&\Z{6}&\text{antipodal pair of} \\
&[0,0,1,\w{3}^2,\w{3}]&\q{20}{12_2}&&\text{eight octahedral} \\
&&&&\text{face-centers on}\\
&&&&\QQ{1}{1},\QQ{1}{2};\\
&&&&\QQ{1}{i} = \LL{2}{5_i} \cap \QQ{}{}\\[\tabGap]

20&[1,1,1,\alpha,\Bar{\alpha}]\ [1,1,1,\Bar{\alpha},\alpha]
&\q{20}{123_1}\ \q{20}{123_2}&\Sym{3}
&\alpha = \frac{-3 +\sqrt{15}\,i}{2}\\[\tabGap]

24&[1,\w{5}^i,\w{5}^j,\w{5}^k,\w{5}^\ell]&\q{24}{ijk\ell}&\Z{5}&
\w{k} = e^{2\,\pi\,i/k}\\[\tabGap]

30&[0,1,i,-1,-i]&\q{30}{1,24_1}=\q{30}{1,35_2}&\Z{4}&
\text{antipodal pair of}\\
&[0,1,-i,-1,i]&\q{30}{1,24_2}=\q{30}{1,35_1}&&
\text{six octahedral}\\
&&&&\text{vertices on}\ \QQ{1}{1} \\[\tabGap]

30&[1,1,\beta,\beta,-2(1+\beta)]
&\q{30}{12,34_1}&\Z{2} \times \Z{2}&\beta = \frac{-2 + \sqrt{5}\,i}{3}\\
&[1,1,\Bar{\beta},\Bar{\beta},-2(1+\Bar{\beta})]&\q{30}{12,34_2}&&\\[\tabGap]

60&[0,1,1,\gamma,\obar{\gamma}]&\q{60}{1,12_1}&\Z{2}&\text{antipodal pair}\\
&[0,1,1,\obar{\gamma},\gamma]&\q{60}{1,12_2}&&\text{of 12 octahedral}\\
&&&&\text{edge-midpoints}\\
&&&&\text{on}\ \QQ{1}{1};\\
&&&&\gamma = -1 + \sqrt{2}\,i

\end{array}$$

\end{table}
%%%%%%%%%%%%%
\begin{table}

\caption{Some fundamental \CP{2} orbits}

\label{tab:planes}

$$\begin{array}{c|c|c|c|c|c|c}

&\text{Algebraic}&\text{Corresponding}&&\text{Set-wise}&\text{Point-wise}&
\text{Restricted}\\
\text{Size}&\text{definition}&\text{point}&\text{Descriptor}&
\text{stabilizer}&\text{stabilizer}&\text{action}\\
\hline

&&&&&&\\[-10pt]

5&\{x_i=0\}&\p{5}{i}&\LL{2}{5_i}&\Sym{4}&\Z{1}&\Sym{4}\\[\tabGap]

10&\{x_i = x_j\}&\p{10}{ij_1}&\LL{2}{10_{ij}}&\Sym{3} \times \Z{2}&\Z{2}&
\Sym{3}\\[\tabGap]

10&\{x_i = -x_j\}&\p{10}{ij_2}&\MM{2}{10_{ij}}&\Sym{3} \times \Z{2}&\Z{1}&
\Sym{3} \times \Z{2}

\end{array}$$

\end{table}
%%%%%%%%%%%%%
\begin{table}

\caption{Special \CP{1} orbits}

\label{tab:lines}

$$\begin{array}{c|c|c|c|c|c}

&&&\text{Set-wise}&\text{Point-wise}&\text{Restricted}\\
\text{Size}&\text{Algebraic definition}&\text{Descriptor}&
\text{stabilizer}&\text{stabilizer}&\text{action}\\
\hline

&&&&&\\[-10pt]

10&\LL{2}{5_i} \cap \LL{2}{5_j}&
\LL{1}{10_{ij}}&\Sym{3} \times \Z{2}&\Z{2}&\Sym{3}\\[\tabGap]

10&\LL{2}{10_{ij}} \cap \LL{2}{10_{jk}} \cap \LL{2}{10_{ik}}&
\MM{1}{10_{ijk}}&\Sym{3} \times \Z{2}&\Sym{3}&\Z{2}\\[\tabGap]

15&\LL{2}{10_{ij}} \cap \LL{2}{10_{k\ell}}\ (i,j \neq k,\ell)
&\LL{1}{15_{ij,k\ell}}
&\D{4}&\Z{2} \times \Z{2}&\Z{2}\\[\tabGap]

15&\MM{2}{10_{ij}} \cap \MM{2}{10_{k\ell}}\ (i,j \neq k,\ell)
&\MM{1}{15_{ij,k\ell}} 
&\D{4}&\Z{2}&\Z{2} \times \Z{2}\\[\tabGap]

30&\LL{2}{5_i} \cap \LL{2}{10_{jk}}\ (i \neq j,k)&\LL{1}{30_{i,jk}}
&\Z{2} \times \Z{2}&\Z{2}&\Z{2}\\[\tabGap]

\end{array}$$

\end{table}
%%%%%%%%%%%
\normalsize

%-------

\section{Parametrized forms} \label{app:forms}
\allowdisplaybreaks

Each case below requires \G{v} invariants to be expressed in terms of
the basic invariants $\Phi_i(v)$.  This amounts to solving a system of
linear equations  whose dimension is that of the respective space of
invariants.  Direct substitution into the basic-invariant expressions
then leads to the descriptions in $K$.

%----------

\subsection*{Basic invariants}

Each $w$-coefficient of $\Phi_\ell(\tau_v w)$ is a degree-$6\,\ell$
invariant in $v$. In terms of $K$, the forms in degrees two and three
are:
\small \begin{align*}
\Phi_{2_K}(w) =&\ \frac{\Phi_2(\tau_v w)}{\Phi(v)^6}\\
=&\ \frac{5}{48}\,\bigl(240\,K_2\,K_3^2\, w_1^2 + 480\,K_1\,K_2\,K_3\,w_1\, w_2 -
48\,K_1^2\, w_2^2 + 240\,K_1^3\,w_2^2 \\
& + 480\,K_1\,K_2\,K_3\,w_1\, w_3 - 96\,K_1\,K_2\,w_2\, w_3 +
480\,K_1\,K_2\,K_3\, w_2\,w_3 - 30\,K_2\,w_3^2 \\
& + 180\,K_1\,K_2\,w_3^2 + 32\,K_2^2\,w_3^2 + 480\,K_2\,K_3^2\,w_1\, w_4 -
60\,K_1\,w_2\,w_4 \\
& + 264\,K_1^2\,w_2\,w_4 + 160\,K_1\,K_2\,w_2\,w_4 - 140\,K_2\,w_3\,w_4 +
184\,K_1\,K_2\,w_3\,w_4 \\
&+ 336\,K_2\,K_3\,w_3\,w_4 - 15\,w_4^2 + 60\,K_1\,w_4^2 + 12\,K_1^2\,w_4^2 +
128\,K_2\,K_3\,w_4^2\bigr) \\
\Phi_{3_K}(w) =&\ \frac{\Phi_3(\tau_v w)}{\Phi(v)^9}\\
=&\ \frac{5}{1728}\,\bigl(-43200\,K_2^2\, K_3^3\,w_1^3 +
25920\,K_1\,K_2\,K_3^2\, w_1^2\,w_2 - 129600\,K_1^2\,K_2\,
K_3^2\,w_1^2\,w_2 \\
& + 51840\,K_1^2\,K_2\,K_3\, w_1\,w_2^2 - 129600\,K_1^2\,K_2\,
K_3^2\,w_1\,w_2^2 + 1944\,K_1^3\,w_2^3 \\
& - 6480\,K_1^4\,w_2^3 - 14400\,K_1^3\,K_2\, w_2^3 +
25920\,K_2^2\,K_3^2\, w_1^2\,w_3 - 129600\,K_2^2\,K_3^3\, w_1^2\,w_3 \\
& + 32400\,K_1\,K_2\,K_3\,w_1\, w_2\,w_3 - 142560\,K_1^2\,K_2\,K_3\,
w_1\,w_2\,w_3  \\
& - 34560\,K_1\,K_2^2\,K_3\, w_1\,w_2\,w_3 + 27432\,K_1^2\,K_2\,
w_2^2\,w_3 - 49680\,K_1^3\,K_2\, w_2^2\,w_3 \\
& - 38880\,K_1^2\,K_2\,K_3\, w_2^2\,w_3 + 37800\,K_2^2\,K_3\,w_1\, w_3^2
- 23760\,K_1\,K_2^2\,K_3\, w_1\,w_3^2 \\
& - 90720\,K_2^2\,K_3^2\, w_1\,w_3^2 + 4860\,K_1\,K_2\,w_2\, w_3^2 -
12960\,K_1^2\,K_2\,w_2\, w_3^2  \\
& - 32400\,K_1^3\,K_2\,w_2\, w_3^2 - 1728\,K_1\,K_2^2\,w_2\, w_3^2 -
17280\,K_1\,K_2^2\,K_3\, w_2\,w_3^2 \\
& + 4860\,K_2^2\,w_3^3 + 3240\,K_1\,K_2^2\,w_3^3 + 384\,K_2^3\,w_3^3 -
9720\,K_2^2\,K_3\, w_3^3 \\
& - 19440\,K_1\,K_2^2\,K_3\, w_3^3 + 16200\,K_2\,K_3^2\, w_1^2\,w_4 -
71280\,K_1\,K_2\,K_3^2\, w_1^2\,w_4  \\
& - 43200\,K_2^2\,K_3^2\, w_1^2\,w_4 + 75600\,K_1\,K_2\,K_3\,w_1\,
w_2\,w_4 - 99360\,K_1^2\,K_2\,K_3\, w_1\,w_2\,w_4 \\
& - 129600\,K_1\,K_2\,K_3^2\, w_1\,w_2\,w_4 + 1620\,K_1^2\,w_2^2\, w_4 -
3888\,K_1^3\, w_2^2\,w_4 \\
& - 6480\,K_1^4\,w_2^2\, w_4 + 17280\,K_1^2\,K_2\, w_2^2\,w_4 -
69120\,K_1^2\,K_2\,K_3\, w_2^2\,w_4 \\
& + 16200\,K_2\,K_3\,w_1\,w_3\, w_4 - 64800\,K_1\,K_2\,K_3\,
w_1\,w_3\,w_4 \\
& - 12960\,K_1^2\,K_2\,K_3\, w_1\,w_3\,w_4 - 86400\,K_2^2\,K_3^2\,
w_1\,w_3\,w_4 + 27000\,K_1\,K_2\,w_2\,w_3\, w_4 \\
& - 48816\,K1^2\,K_2\, w_2\,w_3\,w_4 - 11520\,K_1\,K_2^2\,w_2\,
w_3\,w_4 - 22032\,K_1\,K_2\,K_3\,w_2\, w_3\,w_4 \\
& - 64800\,K_1^2\,K_2\,K_3\, w_2\,w_3\,w_4 + 2025\,K_2\,w_3^2\,w_4 -
3240\,K_1\,K_2\,w_3^2\, w_4 \\
& - 21060\,K_1^2\,K_2\,w_3^2\,w_4 + 2880\,K_2^2\,w_3^2\, w_4 -
7488\,K_1\,K_2^2\, w_3^2\,w_4  \\ 
& - 6912\,K_2^2\,K_3\,w_3^2\, w_4 - 25920\,K_2^2\, K_3^2\,w_3^2\,w_4 +
24300\,K_2\,K_3\,w_1\, w_4^2  \\
& - 48600\,K_1\,K_2\,K_3\,w_1\, w_4^2 - 14400\,K_2^2\,K_3\,w_1\, w_4^2 -
29160\,K_2\,K_3^2\,w_1\, w_4^2 \\
& - 6480\,K_1\,K_2\,K_3^2\, w_1\,w_4^2 + 405\,K_1\,w_2\,w_4^2 -
5508\,K_1^3\,w_2\, w_4^2 \\ 
& + 18000\,K_1\,K_2\,w_2\,w_4^2 - 18720\,K_1^2\,K_2\,w_2\, w_4^2 -
29376\,K_1\,K_2\,K_3\,w_2\, w_4^2  \\ 
& - 25920\,K_1\,K_2\,K_3^2\, w_2\,w_4^2 + 5805\,K_2\,w_3\,w_4^2 -
8640\,K_1\,K_2\,w_3\, w_4^2  \\
& - 3348\,K_1^2\,K_2\,w_3\, w_4^2 - 34992\,K_1\,K_2\,K_3\,w_3\, w_4^2 -
17856\,K_2^2\,K_3\,w_3\, w_4^2 \\
& + 405\,K_1\,w_4^3 - 1620\,K_1^2\,w_4^3 + 324\,K_1^3\,w_4^3 +
3600\,K_2\,w_4^3 - 7200\,K_1\,K_2\,w_4^3 \\
& - 1600\,K_2^2\,w_4^3 - 3456\,K_1\,K_2\,K_3\, w_4^3 -
10368\,K_2\,K_3^2\,w_4^3\bigr)
\end{align*} \normalsize

%----------

\subsection*{Change of coordinates}

Computing the square of the determinant $|\tau_v|$ amounts to
expressing the degree-$20$ invariant $\Psi_{10}(v)^2$ in terms of the
basic forms:
\small \begin{align*}
t_K =& \frac{|\tau_v|^2}{\Phi_2(v)^{24}}
= \frac{\bigl(\Phi_2(v)\,\Phi_3(v)\,\Phi_4(v)\,\Phi_5(v)\,\Psi_{10}(v)\bigr)^2}
{\Phi_2(v)^{24}}\\
=&\ \frac{\Phi_3(v)^4}{\Phi_2(v)^6}\,
\frac{\Phi_4(v)^2}{\Phi_2(v)^4}\,
\frac{\Phi_5(v)^2}{\Phi_2(v)^2\,\Phi_3(v)^2}\,
\frac{\Psi_{10}(v)^2}{\Phi_2(v)^{10}}\\
=&\ \frac{-3125\,K_1^2\,K_2^2\,K_3^2}{13824}\, \bigl( -675 + 9450\,K_1 -
51300\,K_1^2 + 135000\,K_1^3 - 172800\,K_1^4 \\
& + 86400\,K_1^5 + 23700\,K_2 - 147600\,K_1\,K_2 + 111600\,K_1^2\,K_2 +
436800\,K_1^3\,K_2 \\
& - 271800\,K_2^2 + 424800\,K_1\,K_2^2 + 7200\,K_1^2\,K_2^2 +
25600\,K_2^3 - 79200\,K_2\,K_3 \\
& + 535680\,K_1\,K_2\,K_3 - 777600\,K_1^2\,K_2\,K_3 -
576000\,K_1^3\,K_2\,K_3 + 1552320\,K_2^2\,K_3 \\
& - 1238400\,K_1\,K_2^2\,K_3 - 30720\,K_2^3\,K_3 + 68256\,K_2\,K_3^2-
475200\,K_1\,K_2\,K_3^2 \\
& + 864000\,K_1^2\,K_2\,K_3^2 - 3628800\,K_2^2\,K_3^2 +
864000\,K_1\,K_2^2\,K_3^2 \\
& + 4032000\,K_2^2\,K_3^3 - 1728000\,K_2^2\,K_3^4 \bigr).
\end{align*} \normalsize

Each entry of $\tau_v^r\,\tau_v$ is a degree-$12$ invariant in $v$. The
matrix product's expression in $K$ is
\small \begin{align*} 
T_K =&\ \frac{\tau_v^r\,\tau_v}{\Phi_2(v)^6} \\
=&\ \left( \begin{matrix}
240\,K_2\,K_3^2&
2\,K_1\,\bigl( -15 + 66\,K_1 + 40\,K_2 \bigr)\\
200\,K_1\,K_2\,K_3&
48\,K_1\,K_2\,\bigl( -1 + 5\,K_3 \bigr)\\
240\,K_1\,K_2\,K_3&
240\,K_1^2\,\bigl( -1 + 5\,K_1 \bigr)\\
240\,K_2\,K_3^2&
240\,K_1\,K_2\,K_3
\end{matrix} \right.\\
&\left. \begin{matrix}
2\,K_2\,\bigl( -35 + 46\,K_1 + 84\,K_3 \bigr)&
5\,\bigl( -15 + 60\,K_1 + 12\,K_1^2 + 128\,K_2\,K_3 \bigr)\\
2\,K_2\,\bigl( -15 + 90\,K_1 + 16\,K_2 \bigr)&
2\,K_2\,\bigl( -35 + 46\,K_1 + 84\,K_3 \bigr)\\
240\,K_1\,K_2\,\bigl( -1 + 5\,K_3 \bigr)&
2\,K_1\,\bigl( -15 + 66\,K_1 + 40\,K_2 \bigr)\\
240\,K_1\,K_2\,K_3&
240\,K_2\,K_3^2
\end{matrix} \right).
\end{align*} \normalsize
The inverse of $T_K$ results from application of Cramer's rule:
$$T_K^{-1} = \frac{T_K^{\text{cof}}}{|T_K|}$$
where $T_K^{\text{cof}}$ is the matrix of cofactors.

Note that $t_K = |T_K|$.

%----------

\subsection*{Root-selector}

The $w$-coefficients of $\Gamma_v(w)$ are degree-$13$ $v$-invariants.
Expressed in $K$,
\small \begin{align*}
\Gamma_K(w) =&\ \frac{\Gamma_v(w)}{\Phi_2(v)^5\,\Phi_3(v)} \\
=&\ \frac{-125\,\sqrt{5}}{36}\,\bigl(720\,K_2\,K_3^2\, w_1^2 -
288\,K_1\,K_3\,w_1\,w_2 + 1440\,K_1^2\,K_3\,w_1\, w_2 \\
& - 288\,K_1^2\, w_2^2 + 720\,K_1^2\,K_3\,w_2^2 - 288\,K_2\,K_3\,w_1\,w_3 +
1440\,K_2\,K_3^2\,w_1\, w_3\\
& - 180\,K_1\,w_2\,w_3 + 792\,K_1^2\,w_2\,w_3 + 192\,K_1\,K_2\,w_2\,w_3 -
210\,K_2\,w_3^2\\
& + 132\,K_1\,K_2\,w_3^2 + 504\,K_2\,K_3\,w_3^2 - 180\,K_3\,w_1\,w_4 +
792\,K_1\,K_3\,w_1\,w_4\\
& + 480\,K_2\,K_3\,w_1\,w_4 - 420\,K_1\,w_2\,w_4 + 552\,K_1^2\,w_2\,w_4 +
720\,K_1\,K_3\,w_2\,w_4\\
& - 90\,w_3\,w_4 + 360\,K_1\,w_3\,w_4 + 72\,K_1^2\,w_3\,w_4 +
480\,K_2\,K_3\,w_3\,w_4 \\
& - 135\,w_4^2 + 270\,K_1\,w_4^2 + 80\,K_2\,w_4^2 + 162\,K_3\,w_4^2 +
36\,K_1\,K_3\,w_4^2\bigr).
\end{align*} \normalsize

%----------

\subsection*{The $6$-maps}

From the expression for $\phi_6(u)$ in basic invariants and
equivariants, a $K$-parametrized $6$-map $\phi_K(w)$ emerges: 
\small
\begin{align*}
\phi_6(u) =&\ \Phi_2^{15}(v)\,\tau_v\,\Bigl(
2\,\bigl(9\,\Phi_{2_K}(w)\,\Phi_{3_K}(w) -
10\,\Phi_{5_K}(w)\bigr)\,\phi_{1_K}(w)\\ 
&\ - 2\,\bigl(\Phi_{2_K}^2(w) - 5\,\Phi_{4_K}(w)\bigr)\,\phi_{2_K}(w)\\
&\ + 20\,\Phi_{3_K}(w)\,\phi_{3_K}(w) + 15\,\Phi_{2_K}(w)\,\phi_{4_K}(w)
\bigr)\\
=&\ \Phi_2^{15}(v)\,\tau_v\,T_K^{-1}\,\Bigl(
2\,\bigl(9\,\Phi_{2_K}(w)\,\Phi_{3_K}(w) - 10\,\Phi_{5_K}(w)\bigr)\,
\nabla_w^r\Phi_{2_K}(w)\\ 
&\ - 2\,\bigl(\Phi_{2_K}^2(w) - 5\,\Phi_{4_K}(w)\bigr)\,
\nabla_w^r\Phi_{3_K}(w)\\ 
&\ + 20\,\Phi_{3_K}(w)\,\nabla_w^r \Phi_{4_K}(w) +
15\,\Phi_{2_K}(w)\,\nabla_w^r\Phi_{5_K}(w) \Bigr)\\ 
=&\ \Phi_2^{15}(v)\,\tau_v\,\phi_K(w). 
\end{align*} \normalsize

%-------

\section{Basin portraits} \label{app:basins}

\begin{figure}

\resizebox{\basinWidth}{!}{\includegraphics{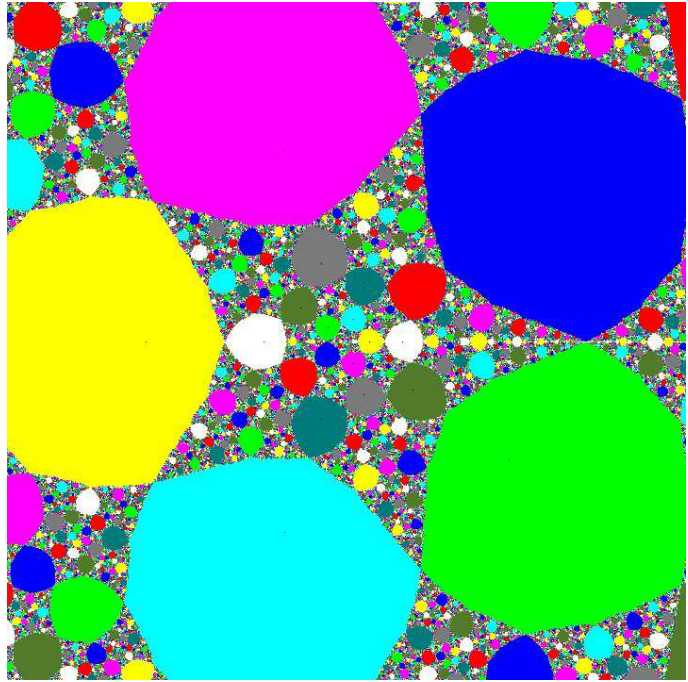}}

\caption{Dynamics of a ruling-preserving 11-map on the quadric's rulings}

\label{fig:dodec11}

\end{figure}

The plots that follow are productions of the program \emph{Dynamics 2}
running on a Dell Dimension XPS with a Pentium II processor. Its BA
process produced \fig{fig:dodec11} while the BAS routine generated the
remaining plots. (See the manual \cite{NY}.) Each procedure
divides the screen into a grid of cells and then colors each cell
according to which attracting point its trajectory approaches. If it
finds no such attractor after $60$ iterates, the cell is black. The BA
algorithm finds the attractor whereas BAS requires the user to specify a
candidate attracting set of points. Each portrait exhibits the highest
resolution available---a $720 \times 720$ grid.  Color versions of the images
appear at \cite{CrassPre}.

In \fig{fig:dodec11}, we have the dodecahedral $11$-map.  Each of the ten pairs 
of antipodal dodecahedral vertices---black
dots---is a period-$2$ superattractor. Their basins fill up
\CP{1} in measure. (Bear in mind that points in the space
of this plot correspond to lines on the quadric surface \QQ{}{}.)

\begin{figure}

\resizebox{\basinWidth}{!}{\includegraphics{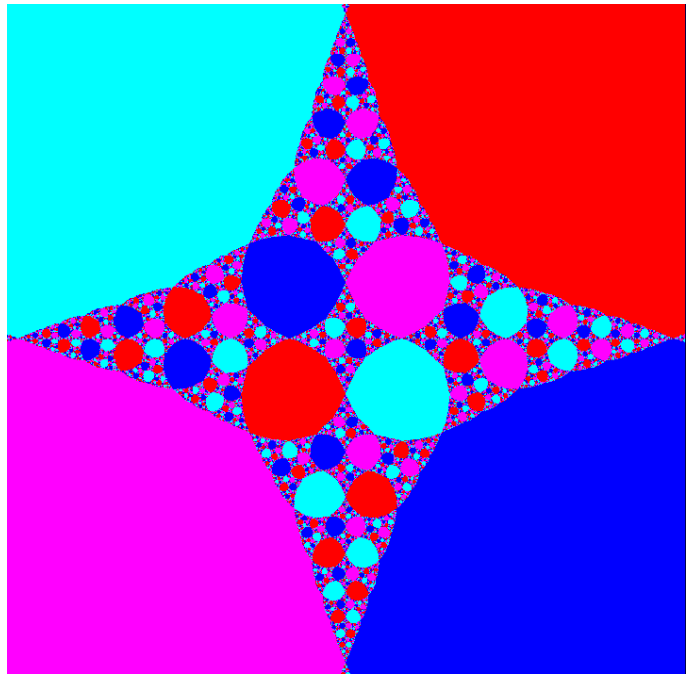}}

\caption{Four basins of attraction for the octahedral $5$-map}

\label{fig:oct5}

\end{figure}

\fig{fig:oct5} indicates the behavior of $h_{11}$ restricted to an
\Sym{4}-symmetric conic \QQ{1}{i}. The four pairs of antipodal
vertices of the cube are period-$2$ superattracting $20$-points whose basins
have full measure on the conic.

\begin{figure}

\resizebox{\basinWidth}{!}{\includegraphics{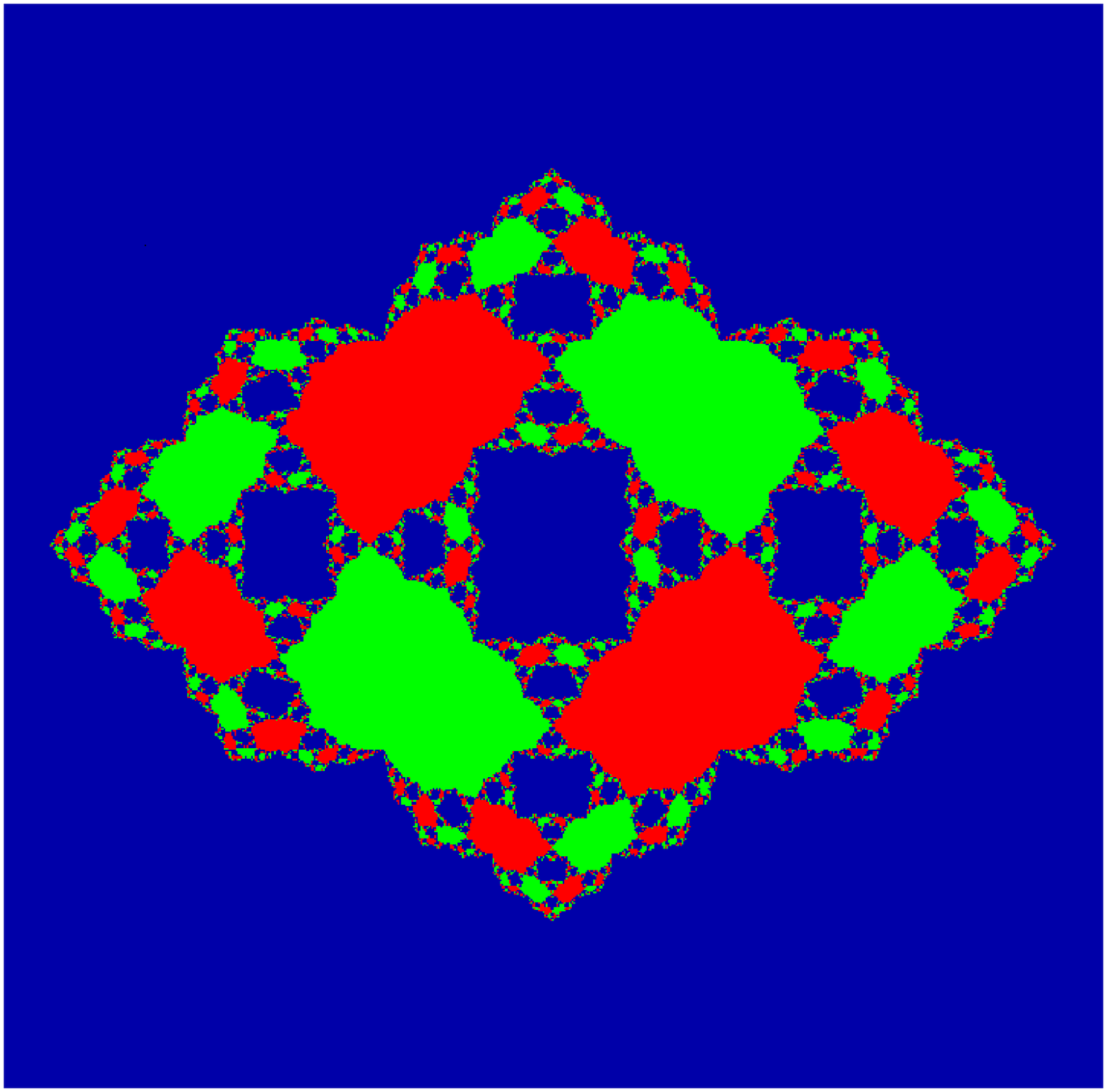}}

\caption{Three basins of attraction for $h_{11}$ restricted to a
$15$-line \LL{1}{15_{ij,k\ell}}}

\label{fig:h11L15}

\end{figure}

\begin{figure}

\resizebox{\basinWidth}{!}{\includegraphics{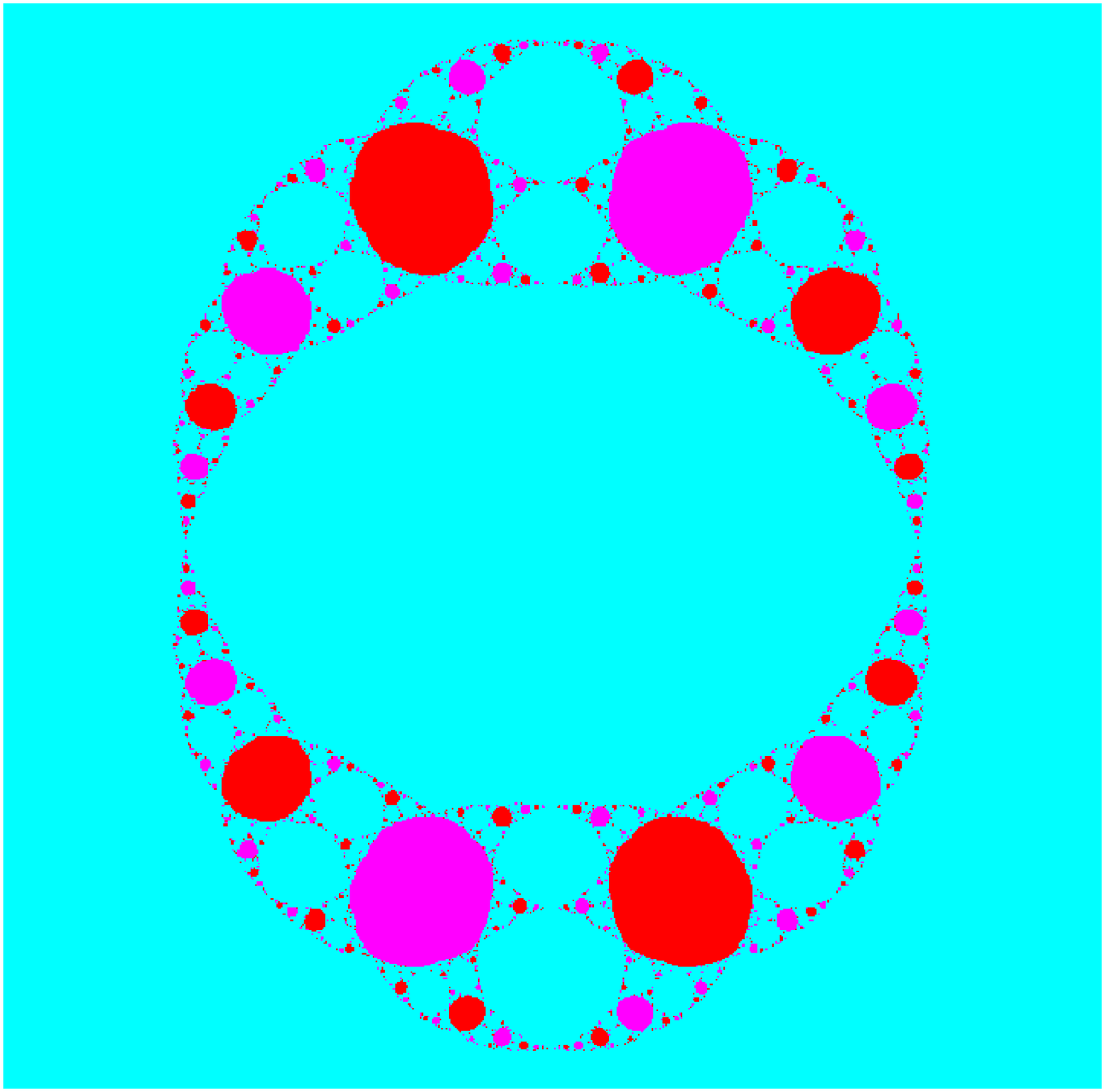}}

\caption{Three basins of attraction for $h_{11}$ restricted to a
$30$-line \LL{1}{30_{i,jk}} }

\label{fig:h11L30}

\end{figure}

\figs~\ref{fig:h11L15} and \ref{fig:h11L30} show the behavior of the octahedral map 
$h_{11}$ on a
$15$-line and a $30$-line respectively.  In the former case, the critical points 
at $0$ and $\infty$ are a pair of $30$-points on
\QQ{}{} that $h_{11}$ exchanges. A pair of fixed $10$-points accounts
for the remaining two basins. At each of these attracting points, the
map repels in at least one direction away from the line. Although the
line has \Z{2} symmetry under \G{120}, the plot displays that of
\mbox{$\Z{2} \times \Z{2}$}. This is a manifestation of an additional
antiholomorphic symmetry
$$x \longrightarrow \obar{x}$$
that extends \G{120} by degree two.

On the $30$-line, the critical points at $0$ and $\infty$ are a pair of octahedral
$60$-points on \QQ{}{} that $h_{11}$ exchanges. The remaining two basins
belong to a pair of $20$-points on \RR{}. At each of these attracting points, the
map repels in at least one direction away from the line. Again,
\mbox{$\Z{2} \times \Z{2}$} symmetry appears.

\begin{figure}

\resizebox{\basinWidth}{!}{\includegraphics{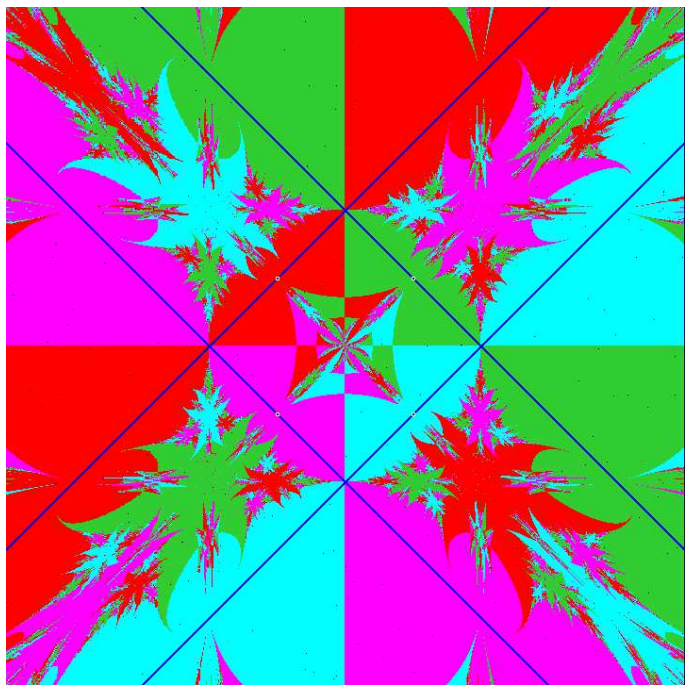}}

\caption{Chaotic attractors for $h_{11}$ on an \RP{2} with \Sym{4}
symmetry}

\label{fig:h11R5}

\end{figure}

\begin{figure}

\resizebox{\basinWidth}{!}{\includegraphics{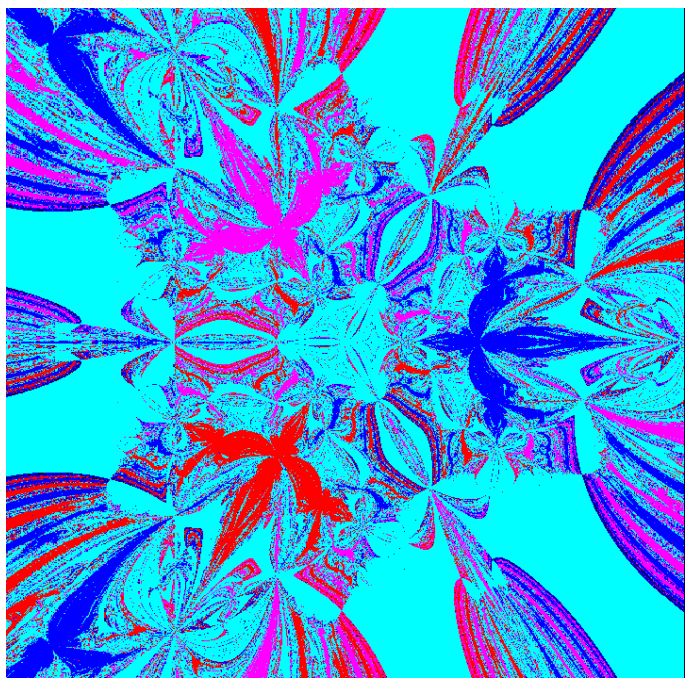}}

\caption{Chaotic attractor for $h_{11}$ on an \RP{2} with \Sym{3}
symmetry}

\label{fig:h11R10}

\end{figure}

In \figs~\ref{fig:h11R5} and \ref{fig:h11R10} we see the restriction of $h_{11}$ 
to an \RP{2} with \Sym{4} symmetry and an \RP{2} with \Sym{3} symmetry.  Each 
case involves a chaotic attractor.  In the former, the attractor consists of the 
four  \RP{1} intersections of \RR,
\LL{2}{5_i}, and the $10$-lines \LL{1}{10_{ij}}. The six intersections
occur at $10$-points \p{10}{k\ell_1} ($k,\ell \neq i$). (In the picture, two of
these intersections occur on the line at infinity.)  The pictured ``lines"
are the images of small circles centered along the edges of the inner
square. This graphical technique specifically relies on the chaotic
behavior of $h_{11}$ along each \RP{1}.

In the \Sym{3} plane, the attracting line is the \RP{1} intersection of \RR, 
\LL{2}{10_{ij}} and the $10$-line
\LL{1}{10_{ij}} at infinity---the light gray basin. The three ``attracting"
$30$-points---they are blowing up---are the vertices 
$$(1,0),\biggl(-\frac{1}{2},\pm\frac{\sqrt{3}}{2}\biggr)$$
of an equilateral triangle about $(0,0)$.

\begin{figure}

\resizebox{\basinWidth}{!}{\includegraphics{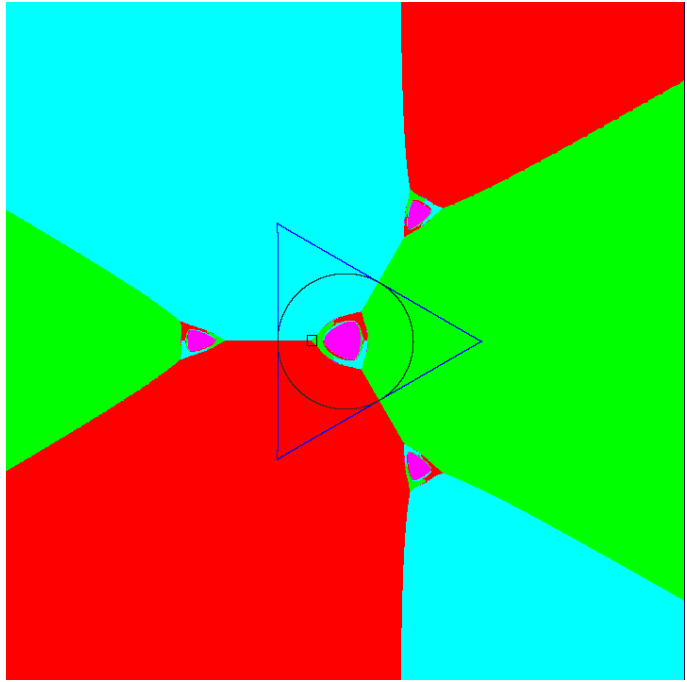}}

\caption{Four basins of attraction for $f_6$ restricted to an \RP{2}}

\label{fig:f6R10}

\end{figure}

\begin{figure}

\resizebox{\basinWidth}{!}{\includegraphics{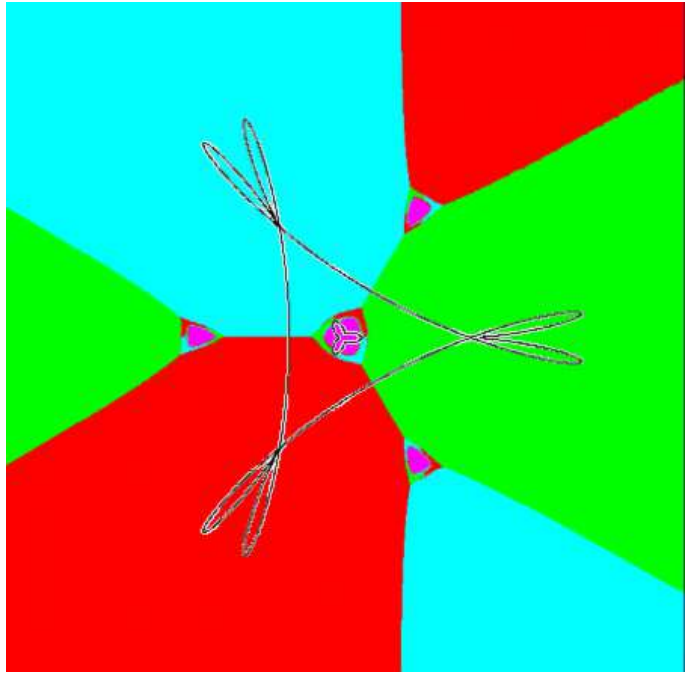}}

\caption{Critical set of $f_6$ restricted to an \RP{2}}

\label{fig:f6R10CSonBasin}

\end{figure}

\begin{figure}

\resizebox{\basinWidth}{!}{\includegraphics{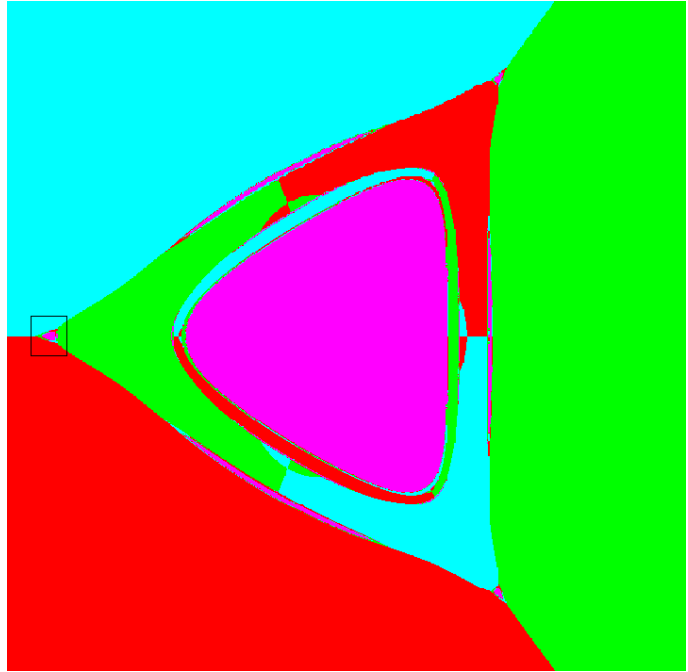}}

\caption{Detail of the left cusp of central basins in \fig{fig:f6R10}}

\label{fig:f6R10mag1}

\end{figure}

\begin{figure}

\resizebox{\basinWidth}{!}{\includegraphics{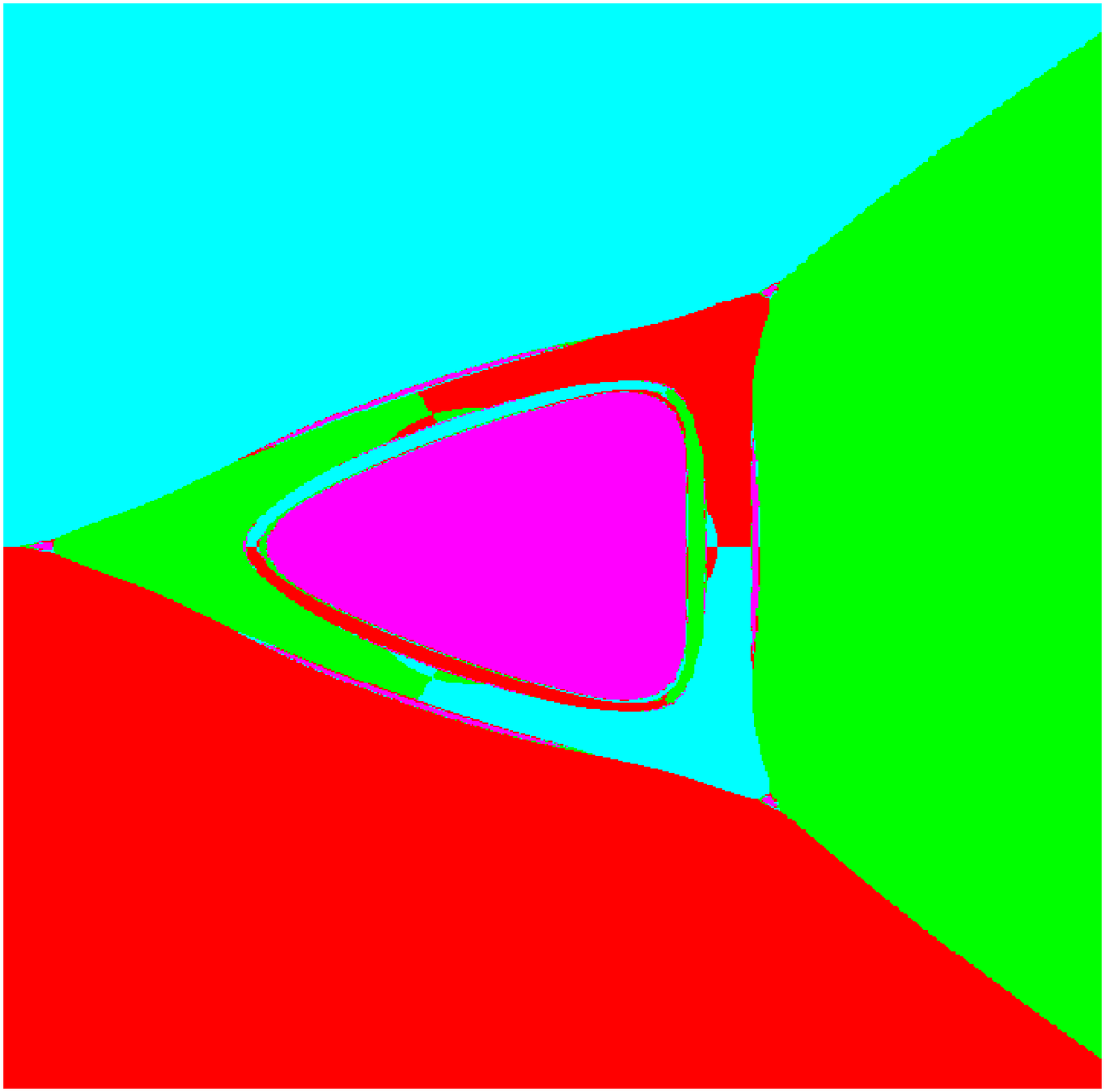}}

\caption{Detail of the left cusp in \fig{fig:f6R10mag1}}

\label{fig:f6R10mag2}

\end{figure}

The remaining images illustrate the dynamics of the quintic-solving $6$-map
$f_6$.  In \figs~\ref{fig:f6R10} through \ref{fig:f6R10mag2}, we see the restriction to
the \RP{2} determined by $\LL{2}{10_{ij}} \cap \RR$.
Since this plane is \Sym{3}-symmetric, the affine
coordinates here are chosen with the three $5$-points at
$$(1,0),\Biggl(-\frac{1}{2},\pm\frac{\sqrt{3}}{2}\Biggr).$$ Three of the
superattracting pipes form a triangle on these points. 
Indeed, the image in \fig{fig:f6R10} of the circle $$\biggl\{x^2 + y^2
=\frac{1}{4}\biggr\}$$ is nearly this triangle. The attractor at $(0,0)$
is the $1$-point orbit \emph{in} the $10$-plane---overall, the
$10$-point \p{10}{ij_2}. In the direction away from the plane, $f_6$
repels at this site along the superattracting pipe \MM{1}{10_{k\ell m}}
($k,\ell,m \neq i,j$). The three ``spokes" at basin boundaries are
pieces of $15$-lines \LL{1}{15_{ij,k\ell}} each of which passes through
a secondary basin that contains a preimage of the central $10$-point.
The boxed region is the approximate content of \fig{fig:f6R10mag1}.

\fig{fig:f6R10CSonBasin} show $h_{11}$'s critical set (\emph{minus} the three ``doubly-critical"
$10$-lines) superimposed on the 
blurry basin portrait. The critical contour is a \emph{Mathematica} plot. Of 
course, the higher order intersections occur at the
$5$-points. All but six critical points appear to belong to the basin of
either a $5$-point or the central $10$-point \p{10}{ij_2}. The six
exceptions lie on the $15$-lines at basin boundaries. If this is so,
then there is no other attracting site---provided that a basin contains 
critical points.

\begin{figure}

\resizebox{\basinWidth}{!}{\includegraphics{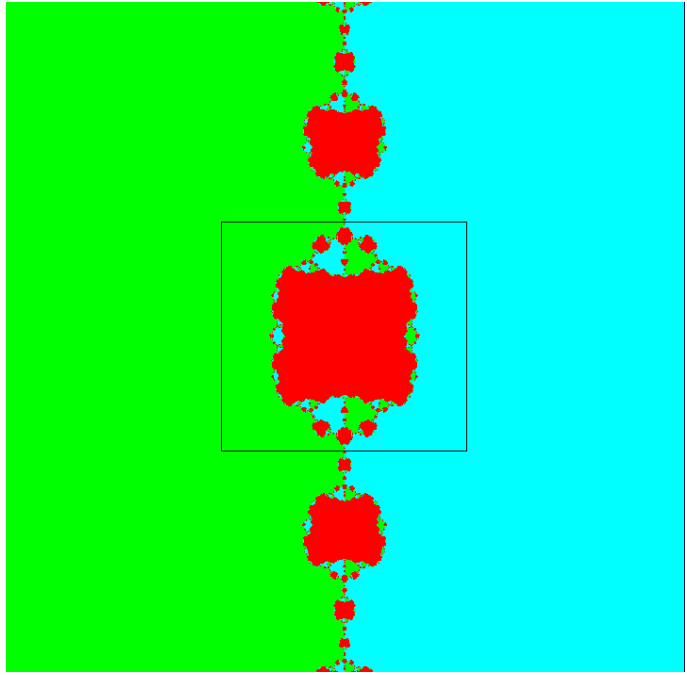}}

\caption{Three basins of attraction for $f_6$ restricted to a $15$-line
\LL{1}{15_{ij,k\ell}}}

\label{fig:f6L15}

\end{figure}

\begin{figure}

\resizebox{\basinWidth}{!}{\includegraphics{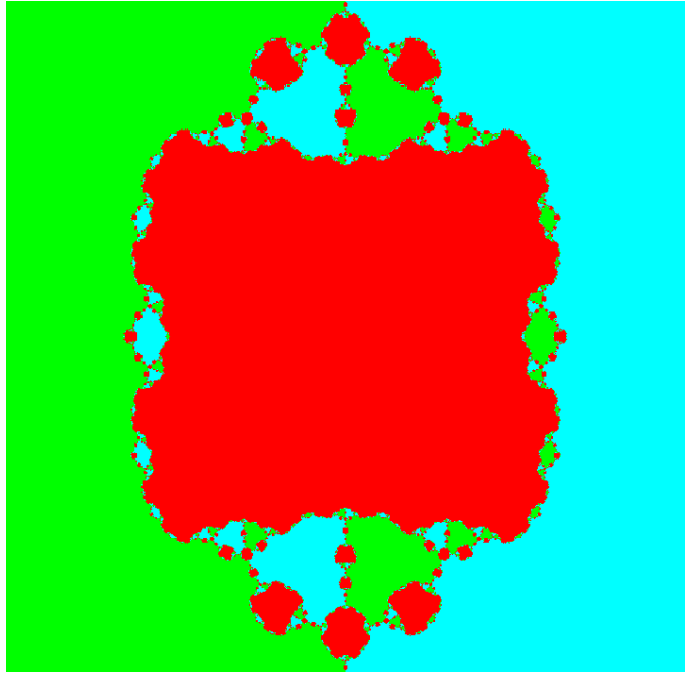}}

\caption{Magnified view of the boxed region in \fig{fig:f6L15}}

\label{fig:f6L15mag}

\end{figure}

In \figs~ \ref{fig:f6L15} and \ref{fig:f6L15mag} we see the map restricted to a $15$-line.
The coordinates of this image place the single $5$-point at $0$ and the
two fixed superattracting $10$-points at $\pm1$. At the latter points,
the map repels in all directions off the line. \fig{fig:f6L15mag}
approximately shows the boxed region.

\begin{figure}

\resizebox{\basinWidth}{!}{\includegraphics{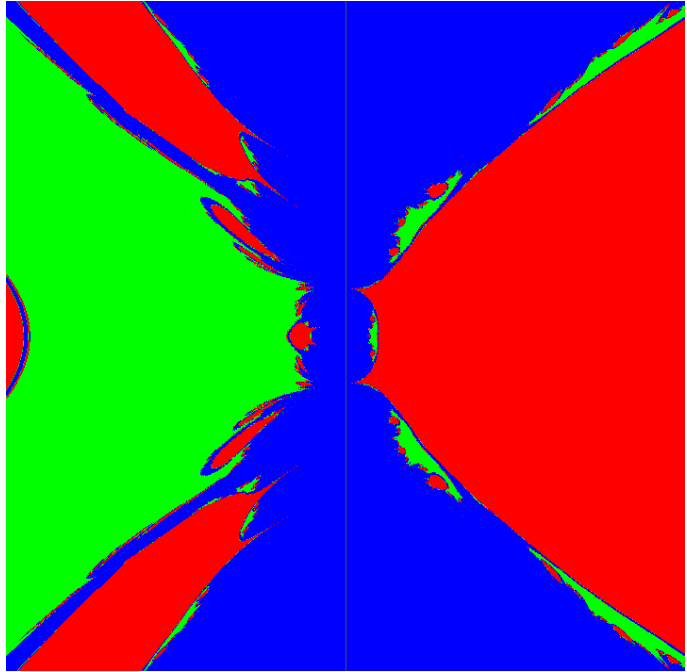}}

\caption{Chaotic attractor for $f_6$ on an \RP{2} with $\Z{2} \times
\Z{2}$ symmetry}

\label{fig:f6R5L10}

\end{figure}

In \fig{fig:f6R5L10}, the space is the \RP{2} intersection of an 
\Sym{4}-invariant
\RP{3} and a $10$-plane \LL{2}{10_{ij}}. The \RP{1} intersection of the
\RP{2} and the $10$-line \LL{1}{10_{ij}} is $\{x=0\}$. By plotting the
trajectory of one of its generic points, this line reveals itself as a
chaotic attractor; the plot shows roughly $20,000$ iterates. The map
attracts at $(1,0),(-1,0)$---the $5$-point \p{5}{k} ($k \neq i,j$) and
$10$-point \p{10}{ij_2} respectively.

%%%%%%%%%%%%%%
%%%%%%%%%%%%%%
%%%%%%%%%%%%%%
%%%%%%%%%%%%%%
%%%%%%%%%%%%%%
%%%%%%%%%%%%%%
%%%%%%%%%%%%%%
%%%%%%%%%%%%%% 
%%%%%%%%%%%%%% 
%%%%%%%%%%%%%%
%%%%%%%%%%%%%%
%%%%%%%%%%%%%%
%%%%%%%%%%%%%%
%%%%%%%%%%%%%%
\clearpage
%-------

%---

\end{document}